\documentclass[twocolumn]{IEEEtran}

\usepackage{amsmath,amssymb,amsthm,epsfig,color,subfigure,empheq,graphicx,graphics,balance,stfloats}
\usepackage{enumerate,url,algorithm,algorithmic,wasysym,epstopdf,enumitem}
\usepackage{multirow,lipsum} 
\usepackage{xcolor}

\usepackage{amsfonts}

\newtheorem{proposition}{Proposition}

\newtheorem{remark}{Remark}

\newcommand\cmb[1]{\textcolor{black}{#1}}

\newcommand \bzero{\mathbf{0}}
\newcommand \bone{\mathbf{1}}

\newcommand \bb{\mathbf{b}}

\newcommand \bd{\mathbf{d}}
\newcommand \be{\mathbf{e}}
\newcommand \bef{\mathbf{f}} 
\newcommand \bg{\mathbf{g}}

\newcommand \bp{\mathbf{p}}
\newcommand \bq{\mathbf{q}}
\newcommand \br{\mathbf{r}}
\newcommand \bs{\mathbf{s}}

\newcommand \bu{\mathbf{u}}
\newcommand \bv{\mathbf{v}}
\newcommand \bw{\mathbf{w}}
\newcommand \bx{\mathbf{x}}

\newcommand \bA{\mathbf{A}}
\newcommand \bB{\mathbf{B}}
\newcommand \bC{\mathbf{C}}
\newcommand \bD{\mathbf{D}}
\newcommand \bE{\mathbf{E}}
\newcommand \bF{\mathbf{F}}
\newcommand \bG{\mathbf{G}}
\newcommand \bH{\mathbf{H}}

\newcommand \bJ{\mathbf{J}}
\newcommand \bK{\mathbf{K}}

\newcommand \bM{\mathbf{M}}
\newcommand \bN{\mathbf{N}}

\newcommand \bP{\mathbf{P}}
\newcommand \bQ{\mathbf{Q}}
\newcommand \bR{\mathbf{R}}

\newcommand \bX{\mathbf{X}}

\newcommand \bgamma{\boldsymbol{\gamma}}

\newcommand \btheta{\boldsymbol{\theta}}

\newcommand \blambda{\boldsymbol{\lambda}}
\newcommand \bmu{\boldsymbol{\mu}}

\newcommand \bphi{\boldsymbol{\phi}}

\newcommand \bGamma{\mathbf{\Gamma}}

\newcommand \bTheta{\mathbf{\Theta}}

\newcommand \mcA{\mathcal{A}}

\newcommand \mcG{\mathcal{G}}

\newcommand \mcL{\mathcal{L}}

\newcommand \mcO{\mathcal{O}}

\newcommand \mcS{\mathcal{S}}

\newcommand \mcU{\mathcal{U}}

\newcommand \mcX{\mathcal{X}}

\newcommand \tbb{\tilde{\mathbf{b}}}

\newcommand \tbx{\tilde{\mathbf{x}}}

\newcommand \tbA{\tilde{\mathbf{A}}}

\newcommand \tbE{\tilde{\mathbf{E}}}

\newcommand \tbN{\tilde{\mathbf{N}}}

\newcommand \tblambda{\tilde{\boldsymbol{\lambda}}}

\newcommand
\tbGamma{\tilde{\mathbf{\Gamma}}}

\newcommand \hbx{\hat{\mathbf{x}}}

\newcommand \hblambda{\hat{\boldsymbol{\lambda}}}

\newcommand \bbb{\bar{\mathbf{b}}}

\newcommand \bbs{\bar{\mathbf{s}}}

\newcommand \bbu{\bar{\mathbf{u}}}

\newcommand \bbA{\bar{\mathbf{A}}}

\newcommand \bbE{\bar{\mathbf{E}}}

\newcommand \bblambda{\bar{\boldsymbol{\lambda}}}

\newcommand \ts{\tilde{s}}

\begin{document}

\title{Fast Probabilistic Hosting Capacity Analysis\\
for Active Distribution Systems}

\author{
	Sina Taheri,~\IEEEmembership{Student Member,~IEEE}, Mana Jalali,~\IEEEmembership{Student Member,~IEEE},\\
Vassilis Kekatos,~\IEEEmembership{Senior Member,~IEEE}, and Lang Tong,~\IEEEmembership{Fellow,~IEEE}
	
%
\vspace*{-2em}
}	
	
\markboth{IEEE TRANSACTIONS ON SMART GRID (submitted February 11, 2020; revised September 18, 2020)}{Taheri, Jalali, Tong, and Kekatos: Fast Probabilistic Hosting Capacity Analysis for Active Distribution Systems}

\maketitle

\begin{abstract}
Interconnection studies for distributed energy resources (DERs) can currently take months since they entail simulating a large number of power flow scenarios. If DERs are to be actively controlled, probabilistic hosting capacity analysis (PHCA) studies become more time-consuming since they require solving multiple optimal power flow (OPF) tasks. PHCA is expedited here by leveraging the powerful tool of multiparametric programming (MPP). Using an approximate grid model, optimal DER setpoints are decided by a quadratic program, which depends on analysis and uncertain parameters in a possibly nonlinear fashion. By reformulating this program, feasible and infeasible OPF instances alike are handled in a unified way to uniquely reveal the location, frequency, and severity of feeder constraint violations. The effect of voltage regulators is also captured by novel approximate models. Upon properly extending MPP to PHCA, we were able to find the exact minimizers for \cmb{518,400 OPF instances on the IEEE 123-bus feeder by solving only 6,905 of them, and 86,400 instances on a 1,160-bus feeder by solving only 2,111 instances.} This accelerated PHCA by a factor of 10. Thus, a utility can promptly infer grid statistics using real-world data without a probabilistic characterization of uncertain parameters.
\end{abstract}


\section{Introduction}
\allowdisplaybreaks
DER integration calls for fast and scalable hosting capacity analysis (HCA) studies. Such studies aim at finding the maximum capacity of DERs on a feeder without violating feeder constraints~\cite{DM17}. Such a deterministic treatment can be overly conservative and does not capture the range of grid conditions. To accommodate the uncertain nature of DERs, probabilistic hosting capacity analysis (PHCA) can be pursued instead. A PHCA study may presume probability distributions on uncertain variables, or rely on real-world data to generate grid scenarios. For any deployment level of DERs, PHCA aims at providing a probabilistic characterization of the attained grid quantities. An example of PHCA is to estimate the cumulative distribution function (cdf) of bus voltages from which risks of violation can be estimated. Since PHCA considers several scenarios to infer reliable grid statistics, it is computationally more demanding than HCA~\cite{DS17}. The computational burden of PHCA may render DER interconnection applications to take months thus hampering the adoption of smart grid technologies. In this context, the objective of this work is to develop fast and scalable solutions for PHCA studies.

P/HCA approaches can be further classified depending on whether DERs operate in a passive or active mode. As an example, solar DERs may operate under maximum power point tracking (MPPT) and perform no or simple (fixed power factor) reactive power control. Under such passive operation, P/HCA studies entail solving a sequence of power flow (PF) tasks. On the other hand, if DER setpoints are actively controlled to minimize a feeder-wide objective, then P/HCA studies involve solving multiple OPF problems.

Commencing with passive DER operation, quasi-static time series (QSTS) analysis is probably the method with the highest accuracy for performing HCA and PHCA alike~\cite{QSTS}. For each DER penetration level, QSTS models time-dependent voltage controllers by running PF solvers on year-long second-based load/solar sequences. QSTS can be expedited by solving the PF tasks at a coarser time granularity and/or by linearizing the PF equations~\cite{DS17}, \cite{HML14}, \cite{NO16}, \cite{TTF12}, \cite{TSPMF18}. The task of maximizing and possibly siting solar generation on a feeder can be also formulated as an optimization problem. Reference \cite{AMZAM18} addresses HCA as a sequence of linear programs by successively linearizing the PF equations. Also, references \cite{DM17} and \cite{Santos17} incorporate regulators, capacitors, and line switches, and handle the related minimizations using a mixed-integer nonlinear solver and a genetic algorithm, respectively. 

Despite the extensive literature on optimal DER operation, P/HCA studies considering DER control are rather limited. To deal with HCA, reference \cite{QWSF17} assumes reactive power control by DERs and maximizes solar penetration subject to upper limits on ohmic power losses and voltage deviations using a linearized grid model. One may alternatively formulate HCA as a mathematical program with equilibrium constraints (MPEC); see \cite{KCR12} for an application of MPEC on strategic investment. Its upper level maximizes solar capacity and the lower level implements DER reactive control subject to voltage constraints. However, the involved mixed-integer program may not scale well to larger feeders. Another way of conducting HCA under active DER operation is based on OPF-controlled DERs. For a specific level of DER penetration, one solves an OPF subject to feeder constraints. If the OPF is (in)feasible, the tested penetration is deemed (in)admissible. The maximum penetration can be readily found by running a sequence of OPFs through bisection. To capture diverse loading and solar conditions, one has to adopt a robust or stochastic OPF~\cite{LDB16}. Such formulations lead to a \emph{single} DER dispatch for all conditions that does not reflect the actual grid operation. Hence, stochastic or robust optimization may be infeasible or its solution too conservative. Another drawback of HCA studies is that they return the maximum capacity, but no sufficient information on the location, frequency, and severity of constraint violations if that maximum capacity is exceeded, such as~\cite{BBDZ18},~\cite{LDB16}

Although PHCA can provide the aforesaid information, it means solving a large number of OPF problems, which are parameterized by uncertain loads and generation. To expedite PHCA, one could leverage the growing literature that uses machine learning to warm-start OPF solvers~\cite{B19}, or infer OPF solutions~\cite{BZ19}. Nonetheless, these works target AC OPF formulations, and thus, cannot take into account the rich problem structure when dealing with OPF tasks that can be posed as convex QPs. Such approximate OPFs are relevant when it comes to planning problems with increased uncertainty.

To accelerate OPF tasks under a linearized grid model, reference~\cite{DM19} trains deep neural networks (DNNs) to predict active constraints, a process that could be improved with the aid of multiparametric programming (MPP); see \cite{BBM03} and \cite{JTT17} for reviews on MPP. Reference \cite{L2OSGC2020} leverages MPP to extract the partial derivatives of OPF solutions with respect to problem parameters and trains a sensitivity-informed DNN that predicts approximate OPF solutions. MPP is also used in \cite{NMRB18} to find the set of active constraints and explores tests for deciding when to stop exploring new MPP solutions rather than actually solving the minimization. MPP has been used in transmission system operations to model congestion and prices in markets~\cite{ZTL11}, \cite{JTT17}; to find the optimal generation investment in electricity markets~\cite{TKV20}; and to trade off load curtailment for reliability in security-constrained economic dispatch~\cite{MBGT19}. To the best of our knowledge, MPP has not been leveraged in distribution system operations.

Previous approaches relying on MPP and active constraint sets like the ones for expediting model predictive control~\cite{BBM03} and~\cite{BBM02}, are tailored to real-time applications and presume relatively few unique constraint patterns. As broader ranges of parameters and larger optimization problems, these \emph{online} approaches are deemed less favorable, especially as faster optimization toolboxes become available. For the PHCA task however, the operator has to deal with a wide range of operating conditions that could yield more constraint patterns. Fortunately however, PHCA is performed \emph{offline} in a batch mode, where the operator can access the complete dataset of OPF instances beforehand. Under this setup, identifying system patterns helps reduce the parameter set with significant overall computational advantages.

The contribution of this work is threefold: \emph{c1)} It develops novel approximate models for voltage regulators to be incorporated in PHCA (Section~\ref{sec:models}); \emph{c2)} Presents a novel formulation of PHCA as a convex quadratic program (QP) to capture both feasible and infeasible instances of the OPF, thus facilitating a quantitative study of feeder constraint violations on a per-bus basis (Section~\ref{sec:problem}); and \emph{c3)} Significantly expedites PHCA by introducing the MPP theory to distribution grid operations (Section~\ref{sec:PP}). Numerical tests using real-world data on a 123- and a 1,160-bus feeder demonstrate that the suggested approach is 10 to 20 times faster than competing alternatives (Section~\ref{sec:tests}).

\emph{Notation:} column vectors (matrices) are denoted by lowercase (uppercase) boldface letters modulo the vector of power flows $\bP+j\bQ$. Calligraphic symbols are reserved for sets. The $n$-th element of $\bx$ is denoted by $x_{n}$; the $(n,m)$-th entry of $\bX$ by $X_{nm}$; and $\|\bx\|_q:=(\sum_{n=1}^N |x_n|^q)^{1/q}$ is the $q$-th norm of $\bx$. Symbols $\bone$, $\bzero$, and $\be_n$ denote the all-ones, all-zeros, and $n$-th canonical vectors. 
\color{black}

\section{Feeder Modeling}\label{sec:models}
\subsection{Modeling Voltage Deviations and Ohmic Losses}\label{subsec:grid}
Consider a single-phase feeder with $N+1$ buses served by the substation bus indexed by $n=0$. Let $v_n$ be the voltage magnitude, and $p_n+jq_n$ the complex power injection at bus $n$. The active power injection $p_n$ is decomposed into $p_n=p_n^g-p_n^c$, where $p_n^g$ is the solar DER generation and $p_n^c$ the load at bus $n$. Reactive injections can be similarly expressed as $q_n=q_n^g-q_n^c$. It is assumed that each bus hosts at most one DER, which captures the aggregation of multiple DERs located on that bus. Collect all but the substation injections and voltages in the $N$-length vectors $\bp=\bp^g-\bp^c$, $\bq=\bq^g-\bq^c$, and $\bv$.

To expedite PHCA, we resort to the approximate grid model of~\cite{BoDo15}, \cite{Deka1}, which is briefly reviewed next. Consider a line connecting buses $m$ and $n$ with impedance $z_{mn}=r_{mn}+jx_{mn}$. If $I_{mn}$ is the line current phasor, the complex voltage drop across the line is $V_m-V_n = z_{mn}I_{mn}$. If $S_{mn} = P_{mn}+jQ_{mn}$ is the complex power flow from buses $m$ to $n$, then $I_{mn}= S_{mn}^*/V_m^*$ by definition. Substituting $I_{mn}$ in the voltage drop equation and multiplying both sides by $V_m^*$ yields
\begin{equation}\label{eq:Vmag}
\left|V_m\right|^2 = V_nV_m^*-S_{mn}^*z_{mn}.
\end{equation}
Maintaining the real part of \eqref{eq:Vmag} and linearizing around the flat voltage profile $\left|V_m^0\right|=\left|V_n^0\right|=V_n^0(V_m^0)^*=1$ provides the approximate voltage drop law
\begin{equation}\label{eq:Vmag-linear}
v_m-v_n \simeq r_{mn}P_{mn}+x_{mn}Q_{mn}.
\end{equation}
Heed the difference between \eqref{eq:Vmag-linear} and the \emph{LinDistFlow} model of~\cite{BW3}, according to which \begin{equation*}
v_m^2-v_n^2 \simeq 2r_{mn}P_{mn}+2x_{mn}Q_{mn}.
\end{equation*}
Note also the \emph{LinDistFlow} model has been derived by ignoring losses and not by linearizing \eqref{eq:Vmag}. Both models have been used widely in grid operations with satisfactory accuracy~\cite{BoDo15}, \cite{VKZG16}. \cmb{We adopt model \eqref{eq:Vmag} rather than \emph{LinDistFlow} so that the operation of voltage regulators is captured by linear inequalities rather than non-convex quadratic constraints in Section~\ref{subsec:regulators}.}

To express voltage drops in a matrix-vector form, stack power flows in vector $\bP+j\bQ$ and partition the branch-bus incidence matrix of the feeder as \cmb{$\tbGamma =[\bgamma_0~\bGamma]$ with $\bGamma^{-1}\bgamma_0=-\bone_N$ since $\tbGamma\bone_{N+1}=\bzero_N$.} If we collect \eqref{eq:Vmag-linear} along all lines and ignore voltage regulators for now, we get
\begin{equation}\label{eq:Av}
\bGamma\bv = \bD_r\bP+\bD_x\bQ-\bgamma_0v_0
\end{equation}
where $v_0$ is the substation voltage and $(\bD_r,\bD_x)$ are diagonal matrices with the values of line resistances and reactances on their main diagonals. Up to a first-order Taylor's expansion around the flat voltage profile, it holds that $S_{mn}\simeq-S_{nm}$. Under this approximation, we get that $\bP=\bGamma^{-\top}\bp$ and $\bQ=\bGamma^{-\top}\bq$. Substitute $\bP$ and $\bQ$ in~\eqref{eq:Av} and pre-multiply by $\bGamma^{-1}$ to finally get~\cite{BoDo16}, \cite{Deka1}
\begin{equation}\label{eq:LDF}
\bv = \bR\bp + \bX\bq + v_0\bone_N
\end{equation}
where $\bR:= \bGamma^{-1}\bD_r\bGamma^{-\top}$ and $\bX:= \bGamma^{-1}\bD_x\bGamma^{-\top}$. Per~\eqref{eq:LDF}, voltages are approximately affine functions of injections. 

We next derive an approximate model for ohmic losses. Under a zero-th or first-order Taylor's series expansion around the flat voltage profile, losses are approximately zero. \cmb{We therefore resort to the following second-order approximation for losses, which is shown in Appendix~\ref{sec:appendixA}.}

\begin{proposition}\label{pro:losses}
Using a second-order Taylor's series expansion at the flat voltage profile, the ohmic losses across a single-phase radial feeder can be approximated as a convex quadratic function of power injections
\begin{equation}\label{eq:LDF-losses}
L\simeq \bp^\top\bR\bp + \bq^\top\bR\bq.
\end{equation}
\end{proposition}

\begin{remark}\label{re:Turitsyn}
Equation \eqref{eq:LDF-losses} has also been used to approximate losses in~\cite{Turitsyn11}. Nonetheless, this formula was obtained in \cite{Turitsyn11} by heuristically approximating $|I_{mn}|^2\simeq P_{mn}^2+Q_{mn}^2$ for every line $(m,n)$, and then approximating $\bP=\bGamma^{-\top}\bp$ and $\bQ=\bGamma^{-\top}\bq$. Proposition~\ref{pro:losses} establishes that this approximation stems in fact from a second-order Taylor's series expansion.
\end{remark}

\cmb{The PHCA methodology proposed in this work could be extended to multiphase feeders. The linearized models of \cite{BoDo15} and \cite{PSCC14} can be used to capture unbalanced conditions, single- and two-phase laterals, and enforce constraints on voltage imbalance. This extension is left for future work. We next model the effect of voltage regulators.}
\color{black}

\subsection{Modeling Voltage Regulators}\label{subsec:regulators}
Modeling a regulator as an ideal transformer for now, its output voltage $v_n$ relates to its input voltage $v_m$ as  $v_n = \alpha_{mn} v_m$. The transformation ratio $\alpha_{mn}$ ranges within $0.9$ and $1.1$ at integer multiples (steps) of $0.0625$; see~\cite{Kersting}. To precisely model these discrete steps, a mixed-integer linear program (MILP) formulation is needed. Given the uncertainty involved in PHCA and for computational simplicity, here we approximate $\alpha_{mn}$ as continuously valued. We distinguish between remotely- and locally-controlled regulators. 

\emph{Remotely-controlled regulators:} For such regulators, the operator can directly select the transformation ratio in near real time, and so $\alpha_{mn}$ can be treated as an optimization variable while dispatching DERs. Under the simplification that $\alpha_{mn}$ varies continuously within $[0.9,1.1]$, there is no need to model its exact value, but it suffices to relate  voltages through
\begin{equation}\label{eq:vmn}
0.9v_m\leq v_n\leq 1.1v_m.
\end{equation}

\begin{figure}[t]
	\centering
	\includegraphics[scale=0.55]{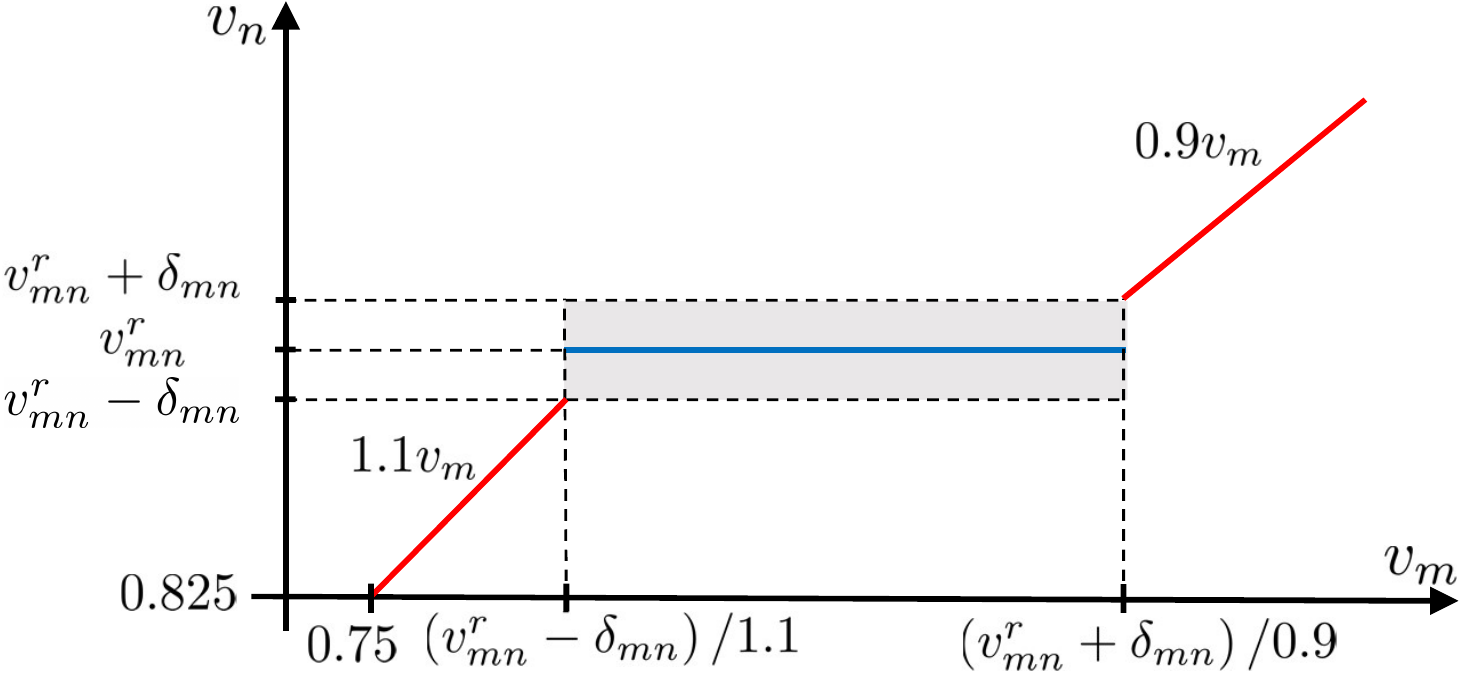}
	\caption{Input-output voltage characteristic for a locally-controlled regulator: The left/rightmost segments occur when taps have reached extreme positions. Within the middle gray box, the output voltage is successfully regulated.}
	\label{fig:dvr}
\end{figure}

\emph{Locally-controlled regulators:} A locally-controlled regulator aims to maintain $v_n$ within a prespecified range $v^r_{mn}\pm \delta_{mn}$ defined by the reference voltage $v^r_{mn}$ and the bandwidth parameter $\delta_{mn}$. The regulator measures its output voltage. If $v_n$ remains outside the range for longer than a given time delay (typically 30-90~sec), the regulator switches its tap position up or down, accordingly. This logic is repeated until $v_n$ is brought within $[v^r_{mn}-\delta_{mn},v^r_{mn}+\delta_{mn}]$, unless an extreme tap position has been reached. Ignoring the time delay, this operation is depicted in Figure~\ref{fig:dvr}. The first segment relates to the case where $v_m$ is quite low and even with the maximum tap position or $\alpha_{mn}=1.1$, the output $v_n=1.1v_m$ remains below $v^r_{mn}-\delta_{mn}$. Similarly, the third segment corresponds to the case where the regulator has attained $\alpha_{mn}=0.9$. Saturation arises due to extreme excursions of the input $v_m$ and should be avoided during normal operation.

Normal operation is captured by the second segment of Figure~\ref{fig:dvr}, where $v_n$ is successfully regulated within $v^r_{mn}\pm \delta_{mn}$. The exact value depends on the tapping sequence and is not known to the operator. For this reason, we propose approximating the shaded area of Fig.~\ref{fig:dvr} with the blue horizontal line passing through the reference voltage, which means
\begin{equation}\label{eq:vn}
v_n\simeq v^r_{mn}.
\end{equation}
Normal operation occurs when the input voltage lies within
\begin{equation}\label{eq:vm}
\frac{v^r_{mn}-\delta_{mn}}{1.1}\leq v_m\leq\frac{v^r_{mn}+\delta_{mn}}{0.9}.
\end{equation} 
If $v_n$ is regulated say within $122\pm 1$~V, then $v^r_{mn}=1.0167$~pu and $\delta_{mn}=0.0083$~pu on a 120-V basis. Hence, the width of the shaded area is $0.2222$~pu, whereas its height is only $0.0167$~pu. This justifies approximating the output voltage by \eqref{eq:vn} when the input voltage remains within \eqref{eq:vm}.

\emph{Regulator with LDC:} Another type of locally-controlled regulators are those equipped with a line drop compensator (LDC)~\cite{Kersting}.
An LDC regulator measures the voltage and current phasors $V_n$ and $I_{mn}$ on its output, and calculates the regulated voltage as
$V_\text{LDC} = V_n-z_\text{LDC}I_{mn}$,
with the LDC impedance setting $z_\text{LDC}=r_\text{LDC}+jx_\text{LDC}$. 
Rather than controlling $v_n=|V_n|$, an LDC regulator controls $v_\text{LDC}=|V_\text{LDC}|$. As in~\eqref{eq:Vmag-linear}, we can approximate $v_\text{LDC} \simeq v_n-r_\text{LDC}P_{mn}-x_\text{LDC}Q_{mn}$ with $P_{mn}+jQ_{mn}$ being the power flow through the regulator. Similar to \eqref{eq:vn}, regulating $v_\text{LDC}$ within $v^r_{mn}\pm\delta_{mn}$ means
\begin{equation}\label{eq:vnLDC}
v_n-r_\text{LDC}P_{mn}-x_\text{LDC}Q_{mn}\simeq v^r_{mn}
\end{equation}
when its input voltage $v_m$ satisfies \eqref{eq:vm}. 

\begin{figure}[t]
	\centering
	\vspace*{-1em}
	\includegraphics[scale=0.3]{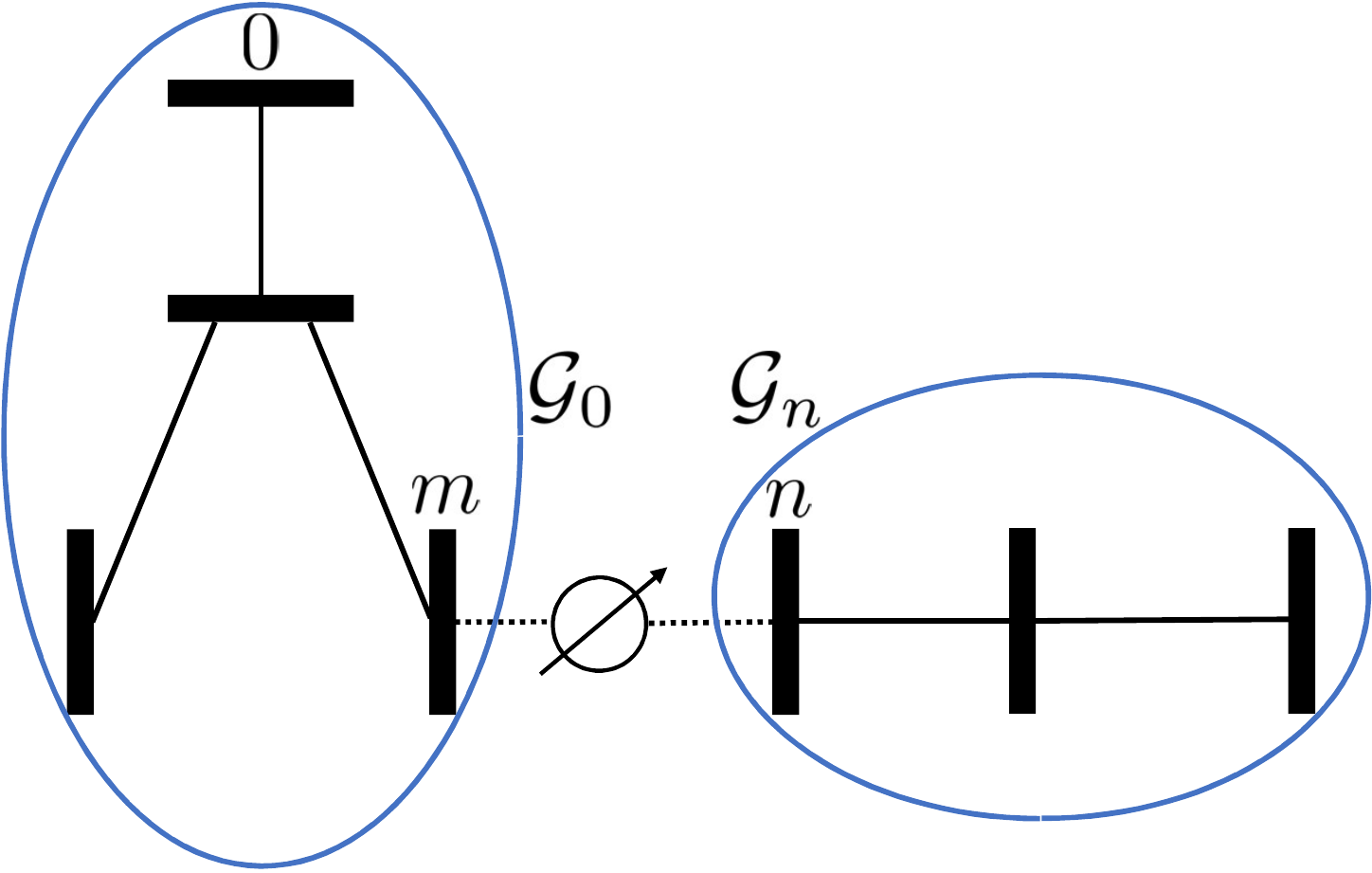}
	\caption{The regulator between buses $m$ and $n$ partitions the grid graph into two subgraphs $\mcG_0$ and $\mcG_n$. Bus $n$ serves as the root for subgraph $\mcG_n$.}
	\vspace*{-1em}
	\label{fig:subgraph}
\end{figure}

Let us now modify the models of Section~\ref{subsec:grid} to account for regulators. Consider the feeder of Figure~\ref{fig:subgraph}. 
The regulator partitions the grid graph into two subgraphs. Within subgraph $\mcG_0$, voltages and losses are modeled by \eqref{eq:LDF}--\eqref{eq:LDF-losses} upon modifying matrices $(\bR,\bX)$ to capture only the edges within $\mcG_0$. Likewise for subgraph $\mcG_n$, with the additional difference that the root voltage $v_0$ in \eqref{eq:LDF} is replaced by $v_n$. The total losses are obviously the sum of \eqref{eq:LDF-losses} over the two subgraphs. 

The two subgraphs are coupled through voltages and power flows. In particular, voltages $v_m$ and $v_n$ are related via \eqref{eq:vmn}--\eqref{eq:vnLDC} depending on the regulator type. The power injection at bus $m$ should be obviously modified as
\begin{equation}\label{eq:pqcouple}
p_m'+jq_m' = (p_m+jq_m)+\sum_{i\in\mcG_n}(p_i+jq_i).
\end{equation}
This partitioning can be straightforwardly generalized to more than one regulators given the tree structure of feeder graphs.

\section{Problem Formulation}\label{sec:problem}
\subsection{Probabilistic Hosting Capacity Analysis (PHCA)}\label{subsec:statement}
A utility considers integrating a specific level of solar DERs on a feeder. The utility may also want to explore the effectiveness of different control options to enable such integration. The computational task of addressing these questions is termed \emph{hosting capacity analysis} and is formally defined next. The operator is assumed to have load estimates, e.g., smart meter data collected every 15 or 60 minutes, and over the course of a year. To keep the analysis tractable, the feeder model stops at the level of distribution transformers, so that a load datum refers to the aggregate load served by a transformer rather than individual households. To model solar generation, the operator collects historical data of solar irradiance over the same year and geographical area and converts them into DER generation.

DERs may also be participating in some form of grid control. We showcase our methodology for the setup where DERs perform reactive power control by following setpoints instructed by the utility. Assuming MPPT, active power generation from solar DERs is not curtailed. To decide the reactive power setpoints for DERs, the utility solves an OPF involving:
\begin{itemize}
\item[\emph{a)}] Optimization variables $\mcO$ (reactive setpoints for DERs, settings for remotely controlled regulators);
\item[\emph{b)}] Uncertain parameters $\mcU$ (loads and solar generation);
\item[\emph{c)}] Analysis parameters $\mcA$ (DER penetration level, settings for locally controlled regulators).
\end{itemize}

For simplicity, we fix $\delta_{mn}$ to the typical value of $\delta_{mn}=0.0083$~pu of $1$V, and leave only $v^r_{mn}$ in $\mcA$. For LDC regulators, the impedance setting $z_\text{LDC}$ can be either fixed or treated as an analysis parameter. The substation voltage $v_0$ is kept in $\mcO$.

Given values for $\mcA$, the operator solves the OPF for all $\mcU$ to obtain $\mcO$. As an example, for $50\%$ penetration $\mcA$, the operator finds the optimal feeder dispatch $\mcO$ for each load/solar scenario experienced over the previous year $\mcU$. For each particular setup of $\mcA$, the operator would like to check whether the OPF is feasible under all scenarios in $\mcU$. If not, what is the probability of an infeasible OPF? For an infeasible OPF instance, what is the minimum allowance in constraint violation one should grant to make it feasible? And if that allowance is granted, which particular constraints got violated and by how much? Our PHCA will answer these questions. We will be able to provide promptly the cdf of voltage violation on any bus and for different solar penetrations. Our PHCA depends on the OPF for deciding $\mcO$, which is formulated next.

\subsection{Optimal Feeder Dispatch }\label{subsec:dispatch}
The set of optimization variables $\mcO$ for our PHCA includes the DER reactive setpoints $\bq^g$ and the substation voltage $v_o$, both collectively denoted by $\bx$. Albeit the ratios $\alpha_{mn}$ for remotely-controlled regulators are optimized, they are not optimization variables. We only enforce \eqref{eq:vmn} and once the OPF is solved, the ratio $\alpha_{mn}$ is obtained as $v_n/v_m$.

Regarding constraints, the reactive power setpoint for DER $n$ is constrained by its apparent power capacity $\overline{s}_n^g$ as
\begin{equation}\label{eq:pv}
|q_n^g|\leq \sqrt{(\overline{s}_n^g)^2 - (p_n^g)^2}
\end{equation}
which is a linear constraint on $\bx$ since $p_n^g$'s are given. Also, power coupling across regulators induces \eqref{eq:pqcouple}. 

We further confine voltage excursions within $\pm 3\%$~pu as
\begin{equation}\label{eq:vlim}
\underline{v}\bone\leq\bv\leq\overline{v}\bone
\end{equation}
with $\underline{v}=0.97$ and $\overline{v}=1.03$. 
The input/output voltages of remotely-controlled regulators should satisfy \eqref{eq:vmn}. The output voltages of locally-controlled regulators are regulated by \eqref{eq:vn} or \eqref{eq:vnLDC}, while their input voltages should satisfy \eqref{eq:vm} to avoid extreme taps. Overall, the feasible set of this OPF is a polytope. 

Given $(\bp^g,\bp^c,\bq^c)\in\mcU$, the DER setpoints $\bq^g$ can be found so they minimize voltage deviations and/or ohmic losses~\cite{jalali19}. Voltage deviations can be captured by the cost function $V(\bx):=\|\bv-\mathbf{1}\|_2^2 + (v_0-1)^2$.
Since voltages $\bv$ are affine functions of $\bx$ from \eqref{eq:LDF}, function $V(\bx)$ is convex quadratic. 
From~\eqref{eq:LDF-losses}, ohmic losses are a convex quadratic function of $\bx$ as well.
The utility may dispatch the feeder by solving
\begin{subequations}\label{eq:opf}
	\begin{align}
	\hbx=\arg\min_{\bx} \ &~F(\bx):=\beta V(\bx) + (1-\beta) L(\bx)\label{eq:opf:cost}\\
	\mathrm{s. to} \ &~\bC_1\bx\leq\bd_1\label{eq:opf:s1}\\
	&~\bC_2\bx\leq\bd_2.\label{eq:opf:s2}
	\end{align} 
\end{subequations}
where $\beta\in[0,1]$ is decided by the utility to trade voltage regulation for power losses. Problem~\eqref{eq:opf} is a linearly-constrained convex quadratic program. For reasons explained in Section~\ref{subsec:feasible}, we have grouped constraints into two subsets:

\emph{Subset \eqref{eq:opf:s1}} consists of linear (in)equalities that encode the actual feeder operation and cannot be relaxed, namely \eqref{eq:pv} and the regulator-related constraints \eqref{eq:vmn}, \eqref{eq:vn}, \eqref{eq:vnLDC}, and \eqref{eq:pqcouple}.

\emph{Subset \eqref{eq:opf:s2}} includes linear inequalities encoding feeder constraints the DER dispatch should satisfy, such as \eqref{eq:vm} and \eqref{eq:vlim}. Such constraints can be relaxed to make the dispatch feasible. It is exactly these constraints, the utility would like to study and obtain the cdf's of their violations during PHCA.

It is worth pointing out that assigning a constraint to \eqref{eq:opf:s1} or \eqref{eq:opf:s2} may depend on the goals of a particular PHCA. For instance, the apparent constraint of \eqref{eq:pv} stems from overheating/lifetime concerns and may be relaxed for short periods of time. Hence, the operator may consider moving \eqref{eq:pv} from \eqref{eq:opf:s1} to \eqref{eq:opf:s2}. Vice versa, to avoid saturation and correctly capture regulator's operation, the utility may decide transferring constraint \eqref{eq:vm} from \eqref{eq:opf:s2} to \eqref{eq:opf:s1}. 

The utility would like to solve \eqref{eq:opf} over multiple instances of $\mcU\times \mcA$. The analysis parameters may involve the settings $(v_{mn}^r,z_\text{LDC})$ for locally-controlled regulators, parameter $\beta$, or the level of solar penetration. The reason for partitioning its feasible set is to handle jointly (in)feasible instances of \eqref{eq:opf} over $\mcU\times \mcA$, as pursued next. 

\subsection{Jointly Handling Feasible and Infeasible OPF Instances}\label{subsec:feasible}
Under some instances in $\mcU\times \mcA$, problem~\eqref{eq:opf} may become infeasible. One approach to pinpoint the cause of infeasibility is to allow for voltage violations by introducing a slack optimization variable into \eqref{eq:opf:s2} and penalize its effect as
\begin{subequations}\label{eq:opfr}
\begin{align}
\left(\tbx,\ts\right)=\arg\min_{\bx,s\geq 0}~&~F(\bx) +\nu s^2 + \eta s\label{eq:main:a}\\
\mathrm{s.to} ~&~\bC_1\bx\leq\bd_1\label{eq:opfr:s1} \\
~&~\bC_2\bx\leq\bd_2+s\bone.\label{eq:opfr:s2}
\end{align}
\end{subequations}
Voltage violations can occur on lower or upper limits, and in one or multiple buses. The slack variable $s$ is penalized by a linear-quadratic cost determined by positive parameters $(\eta,\nu)$. \cmb{Notice that \eqref{eq:opfr} is always feasible with $\bq^g=\bzero$.}

To explore the connection between \eqref{eq:opf} and \eqref{eq:opfr}, consider an OPF instance in $\mcU\times \mcA$ and let $\hbx$ and $(\tbx,\ts)$ be the minimizers of \eqref{eq:opf} and \eqref{eq:opfr}, respectively. For sufficiently large $(\eta,\nu)$, the minimizer of \eqref{eq:opfr} exhibits three properties:

\emph{p1)} If \eqref{eq:opf} is feasible, then \eqref{eq:opfr} yields $\tbx=\hbx$ and $\ts=0$.
	
\emph{p2)} If \eqref{eq:opf} is infeasible, one may try loosening the constraints in \eqref{eq:opf:s2} by $\ts$. The solution to this relaxation of \eqref{eq:opf} would be $\hbx=\tbx$. This is easy to verify by fixing $s=\ts$ in \eqref{eq:opfr}, and then minimizing \eqref{eq:opfr} over $\bx$ to get $\tbx$. Therefore, for infeasible instances of \eqref{eq:opf}, problem \eqref{eq:opfr} yields a relaxed feeder dispatch that if actually implemented, the maximum constraint violation would be equal to or smaller than $\ts$. 

\emph{p3)} If \eqref{eq:opf} is infeasible, scalar $\ts$ is the minimum value by which we have to relax \eqref{eq:opf} to make it feasible.

When solving \eqref{eq:opfr} with an actual solver however, setting $(\eta,\nu)$ at large values could jeopardize numerical stability. Hence, one would prefer the smallest values of $(\eta,\nu)$ that still achieve properties \emph{p1)}--\emph{p3)}. First note that \emph{p2)} holds for all positive $\eta$ and $\nu$. Proceeding with \emph{p1)}, the next proposition (shown in Appendix~\ref{sec:appendixB}) provides a lower bound on $\eta$ to ensure that \emph{p1)} still holds for any feasible instance of \eqref{eq:opf}.

\color{black}
\begin{proposition}\label{prop:exact penalty}
	Consider the optimization problem 
		\begin{equation}\label{eq:P1}
		\hbx = \arg\min_{\bx}\left\{f(\bx):g_i(\bx)\leq0,~ i=1,\dots,K\right\}\tag{P1}
		\end{equation}		
where $f$ is a strictly convex function and $g_i$'s are convex functions of $\bx$. Assume Slater's condition holds for \eqref{eq:P1}, and let $\hblambda$ be the vector of optimal Lagrange multipliers. Consider also problem 
		\begin{align}\label{eq:P2}
		\left(\tbx,\ts\right) =\arg\min_{\bx,s\geq0}~&~f(\bx)+p(s)\tag{P2}\\
		\mathrm{s. to}~&~g_i(\bx)\leq s,\quad i=1,\dots,K\nonumber
		\end{align}
where $p(s)$ is a differentiable and increasing convex function. If $\frac{d p(s)}{ds}>\hblambda^\top\bone$ for all $s\geq 0$, then $\tbx=\hbx$ and $\ts=0$.
\end{proposition}

To apply Proposition~\ref{prop:exact penalty} to \eqref{eq:opfr}, observe that the derivative of $p(s)=\nu s^2 +\eta s$ with respect to $s$ is $2\nu s + \eta\geq \eta$ for all $s\geq 0$. Let $\hblambda$ be the optimal Lagrange multipliers for constraint \eqref{eq:opf:s2} for a feasible instance of \eqref{eq:opf}. According to Proposition~\ref{prop:exact penalty}, property \emph{p1)} holds if we set $\eta\geq \bone^\top\hblambda$. Two remarks on the penalization cost of \eqref{eq:opfr} are now in order.

First, Proposition~\ref{prop:exact penalty} predicates that a purely linear penalty would suffice. The motivation behind our linear-quadratic penalty is to make \eqref{eq:opfr} a strictly convex QP, so it is readily amenable to MPP later. Given \eqref{eq:opf} is an QP already, appending linear-quadratic penalties does not change the problem class. 

Second, to keep our PHCA tractable, we would like to keep $(\nu,\eta)$ constant over all problem instances. To use Proposition~\ref{prop:exact penalty} though, we need to find the maximum Lagrange multiplier sum over all problem instances. Since this may be hard, we resorted to solving a few instances of \eqref{eq:opf}; recording the maximum sum; and increasing that by an order of magnitude to set $\eta$. 

Regarding \emph{p3)}, we have no analytical claim. We numerically observed that \emph{p3)} is attained if $\nu$ is chosen larger or equal to the largest eigenvalue of the Hessian matrix of $F(\bx)$. Since solving \eqref{eq:opfr} for thousands of instances can be time consuming, we next leverage its structure and devise fast PHCA solvers.

\section{Fast Data-Based Feeder Analysis}\label{sec:PP}

\subsection{Multiparametric Programming (MPP)}\label{subsec:MPP}
MPP studies parameterized optimization problems; see \cite{TJB03} and \cite{JTT17}. Given a parameter vector $\btheta\in\bTheta\subseteq \mathbb{R}^K$, consider the general quadratic program over variable $\bx\in\mathbb{R}^N$
\begin{subequations}\label{eq:qp}
	\begin{align}
	\underset{\bx}{\min}~&~\frac{1}{2}\bx^\top\bH\bx+\left(\bC\btheta+\bd\right)^\top\bx\label{eq:qp:cost}\\
	\mathrm{s.to}~&~\bA\bx\leq\bE\btheta+\bb&&:\blambda\label{eq:qp:ineq}\\
			      &~\bB\bx=\bF\btheta+\bef&&:\bmu\label{eq:qp:eq}
	\end{align}
\end{subequations}
where $\bH\succeq\bzero$; $\mathbf{A}\in\mathbb{R}^{I\times N}$; $\mathbf{B}\in\mathbb{R}^{E\times N}$; and the rest of matrices/vectors are of conformable dimensions. Vectors $\blambda$ and $\bmu$ collect the dual variables of \eqref{eq:qp}. If $\bH=\bzero$, problem \eqref{eq:qp} is a multiparametric linear program (MPLP). Otherwise, it is a multiparametric quadratic program (MPQP). Either way, the subset of $\bTheta$ for which \eqref{eq:qp} is feasible, can be partitioned into distinct regions, termed \emph{critical regions}~\cite{BBM03}. Interestingly enough, each region is described by a polytope in $\bTheta$, and within each region the primal/dual solutions to \eqref{eq:qp} can be expressed as affine functions of $\btheta$. 

\cmb{Different from the typical MPQP setup where parameters appear only on the right-hand side of linear inequalities~\cite{BBM02},~\cite{BBM03}, problem \eqref{eq:qp} involves $\btheta$ in its objective and has the linear equalities of \eqref{eq:qp:eq} as well. To describe the solution of~\eqref{eq:qp}, we next review MPP results for QPs with parameters in the cost~\cite{EMPC}. We start with its Lagrangian function}
\begin{align}
\mcL(\bx;\blambda,\bmu) &= \frac{1}{2}\bx^\top\bH\bx+\left(\bC\btheta+\bd\right)^\top\bx+\blambda^\top(\bA\bx-\bE\btheta-\bb)\nonumber\\
&+\bmu^\top\left(\bB\bx-\bF\btheta-\bef\right).\nonumber
\end{align} 

Suppose \eqref{eq:qp} is solved for a particular $\btheta\in\bTheta$. Let matrix $\tbA$ be obtained by sampling the rows of $\bA$ corresponding to the constraints in~\eqref{eq:qp:ineq} satisfied with equality at the optimum, i.e., the \emph{active constraints}. Matrix $\tbE$ and vectors $\tbb$ and $\tblambda$ are defined similarly. Let matrix $\bbA$ collect the remaining rows of $\bA$, that is the rows corresponding to the \emph{inactive} constraints in~\eqref{eq:qp:ineq} satisfied with strict inequality. Define $(\bbE,\bbb,\bblambda)$ similarly. The optimal primal/dual variables should satisfy the Karush-Kuhn-Tucker (KKT) conditions
\begin{subequations}\label{eq:cond}
	\begin{align}
		&\bH\bx+\bC\btheta+\bd+\tbA^\top\tblambda +\bB^\top\bmu= \bzero\label{eq:cond:lo}\\
		&\tbA\bx=\tbE\btheta+\tbb\label{eq:cond:pf1}\\
		&\bB\bx = \bF\btheta+\bef\label{eq:cond:pf2}\\
		&\bbA\bx<\bbE\btheta+\bbb\label{eq:cond:pf3}\\
		&\tblambda\geq\bzero,\label{eq:cond:df}\\
		&\bblambda=\bzero.\label{eq:cond:cs}
	\end{align}
\end{subequations}
If $\bH\succ\bzero$, the primal minimizer is expressed from \eqref{eq:cond:lo} as
\begin{equation}\label{eq:primal_sol}
\bx = -\bH^{-1}(\bC\btheta+\bd+\tbA^\top\tblambda+\bB^\top\bmu).
\end{equation}
Define matrix $\bK:=[\tbA^\top~~\bB^\top]^\top$. If $\bK$ is full row-rank, then $\bK\bH^{-1}\bK^\top\succ \bzero$. Substituting~\eqref{eq:primal_sol} into~\eqref{eq:cond:pf1}--\eqref{eq:cond:pf2}, the optimal dual variables are expressed as
\begin{align}
\begin{bmatrix}
\tblambda\\
\bmu
\end{bmatrix} =\begin{bmatrix}
\bG_1\\
\bG_2
\end{bmatrix}\btheta+\begin{bmatrix}
\bw_1\\
\bw_2
\end{bmatrix}\label{eq:dual_sol}
\end{align}
where
\begin{align}
&\begin{bmatrix}
\bG_1\\\bG_2\end{bmatrix}:=-\left(\bK\bH^{-1}\bK^\top\right)^{-1}\left(\bK\bH^{-1}\bC+\begin{bmatrix}\tbE\\\bF\end{bmatrix}\right)\label{eq:G1G2}\\
&\begin{bmatrix}
\bw_1\\\bw_2\end{bmatrix}:=-\left(\bK\bH^{-1}\bK^\top\right)^{-1}\left(\bK\bH^{-1}\bd+\begin{bmatrix}\tbb\\\bef\end{bmatrix}\right)\label{eq:w1w2}.
\end{align}
Substituting~\eqref{eq:dual_sol} in \eqref{eq:primal_sol} yields the optimal primal variable
\begin{align}\label{eq:primal-sol2}
\bx =\bM\btheta+\br
\end{align}
where $\bM := -\bH^{-1}(\bC+\tbA^\top\bG_1+\bB^\top\bG_2)$ and $\br := -\bH^{-1}(\bd+\tbA^\top\bw_1+\bB^\top\bw_2)$.

Clearly from \eqref{eq:dual_sol}--\eqref{eq:primal-sol2}, the primal and dual solutions to \eqref{eq:qp} are affine functions of $\btheta$. Interestingly, the same affine functions apply for any other $\btheta$ yielding the same active constraints. Conversely, given a set of constraint indexes, conditions \eqref{eq:cond:pf3}--\eqref{eq:cond:df} can be used to obtain the subset of $\btheta$'s that make those indexes active~\cite{TJB03}. Let sets $\mcS_p$ and $\mcS_d$ comprise respectively all $\btheta$'s satisfying \eqref{eq:cond:pf3} and \eqref{eq:cond:df} for a given index set of active constraints. Using \eqref{eq:dual_sol}--\eqref{eq:primal-sol2}, the sets $\mcS_p$ and $\mcS_d$ can be expressed as polyhedra in $\bTheta$
\begin{subequations}\label{eq:cr}
\begin{align}
&\mcS_p := \left\lbrace\btheta|\left(\bbA\bM-\bbE\right)\btheta<\bbb-\bbA\br\right\rbrace\label{eq:cr:p}\\
&\mcS_d:= \left\lbrace\btheta|\bG_1\btheta+\bw_1\geq\bzero\right\rbrace.\label{eq:cr:d}
\end{align}
\end{subequations}
The intersection $\mcS_p\cap\mcS_d\subseteq \bTheta$ defines a \emph{critical region}, that is a subset of $\btheta$'s yielding the same active constraints for~\eqref{eq:qp}. Each critical region can be uniquely characterized by a set of indexes of active constraints. Inside each critical region, the optimal primal solution is an affine function of $\btheta$, such as the one in \eqref{eq:primal-sol2} with $(\bM,\br)$ changing per region. 

The key feature of MPP is that once a critical region has been visited, solving \eqref{eq:qp} for \emph{any} new $\btheta$ within that region becomes trivial: The primal and dual solutions to \eqref{eq:qp} for the new $\btheta$ are readily provided by \eqref{eq:dual_sol}--\eqref{eq:primal-sol2}. An MPQP can have in theory exponentially many critical regions. Nonetheless, the MPQPs involved in practical applications (e.g., model predictive control~\cite{TJB03}; or power transmission system operations~\cite{ZTL11}, \cite{JTT17}, \cite{MBGT19}), oftentimes exhibit a limited number of regions. We next leverage MPP to expedite PHCA.

\subsection{PHCA as a Multiparametric Quadratic Program (MPQP)}
The PHCA task of~\eqref{eq:opfr} can be interpreted as an MPQP of the form in~\eqref{eq:qp} over $\bx=[\bq^\top~v_0~s]^\top$. Its parameter vector $\btheta$ is constructed from instances of $\mcU\times \mcA$. Suppose a sample $u_s=(\bp^c,\bq^c,\bp^g)\in\mcU$, where vector $\bp^g$ is in reference to $100\%$ solar penetration. Suppose also that $\mcA$ contains only the solar penetration level. If one wants to solve \eqref{eq:opfr} for the uncertain parameters $u_s$ but under $50\%$ penetration, the corresponding parameter vector $\btheta_s$ includes $(\bp^c,\bq^c,0.5\cdot\bp^g)$. It further includes the quantities $\sqrt{(\overline{s}_n^g)^2 - (0.5\cdot p_n^g)^2}$ for all $n$, since they appear in \eqref{eq:pv}. Even though $\btheta$ appears linearly in \eqref{eq:qp}, it can be a nonlinear mapping from $\mcU\times \mcA$ to $\bTheta$.

Having defined $(\bx,\btheta)$, the matrices/vectors of \eqref{eq:qp} can be obtained from \eqref{eq:opfr}. By including the penalty $\nu s^2$ in~\eqref{eq:opfr}, we ensure the related $\bH$ is positive definite. To maintain the condition number of $\bH$ within a reasonable range, $\nu$ was set equal to the largest eigenvalue of $\bR$. By definition of $F(\bx)$, it is not hard to verify the quadratic component of the cost in \eqref{eq:opf} is independent of $\btheta$, whereas its linear one depends affinely on $\btheta$ thus complying with the parametric form in \eqref{eq:qp:cost}.

\begin{algorithm}[t]
\caption{Multiparametric OPF (MP-OPF)}\label{alg:MP-OPF}
\begin{algorithmic}[1]
\renewcommand{\algorithmicrequire}{\textbf{Input:}}
\renewcommand{\algorithmicensure}{\textbf{Output:}} 
\REQUIRE Set of OPF scenarios $\bTheta=\left\{\btheta_s\right\}_{s=1}^S$
\ENSURE OPF solutions $\{\bx_s\}_{s=1}^S$ to \eqref{eq:opfr} via \eqref{eq:qp} for $\btheta=\btheta_s$
		\WHILE {$\bTheta\neq \emptyset$}
		\STATE Randomly select $\btheta_o\in\bTheta$ and $\bTheta\gets\bTheta\setminus\btheta_o$
		\STATE Solve~\eqref{eq:qp} for $\btheta=\btheta_o$ to find the optimal $\bx_o$ and its active constraints
		\IF {matrix $\bK$ is full row-rank,}
			\STATE Compute region parameters $(\bG_1,\bG_2,\bw_1,\bw_2,\bM,\br)$
			\STATE Compute polytope $\mcS_p\cap\mcS_d$ describing this region
			\FOR {all $\btheta_s\in\bTheta$}
				\IF {$\btheta_s\in\mcS_p\cap\mcS_d$,}
					\STATE Find OPF solution as $\bx_s=\bM\btheta_s+\br$
					\STATE $\bTheta\gets\bTheta\setminus\btheta_s$ 
				\ENDIF
			\ENDFOR
		\ENDIF
		\ENDWHILE
	\end{algorithmic}
\end{algorithm}

The operator would like to solve \eqref{eq:opfr} for a large number of OPF scenarios comprising set $\bTheta:=\{\btheta_s\}_{s=1}^S$ with related OPF minimizers $\{\bx_s\}_{s=1}^S$. Instead of solving \eqref{eq:opfr} for each $\btheta_s\in\bTheta$, we can utilize MPP and solve \eqref{eq:opfr} only for one $\btheta_s$ per critical region of $\bTheta$. The process is tabulated as Algorithm~\ref{alg:MP-OPF} and is henceforth termed \emph{MP-OPF}. At Step 3, MP-OPF solves \eqref{eq:opfr} via \eqref{eq:qp} for a particular $\btheta_o\in\bTheta$. For the visited critical region, MP-OPF computes its parameters $(\bM,\br)$ and polytopic description $\mcS_p\cap \mcS_d$. To avoid storing $(\bM,\br)$ and $\mcS_p\cap \mcS_d$ for all regions, MP-OPF proceeds by scanning through $\bTheta$ to identify other $\btheta$'s belonging to the just visited region. For those $\btheta_s$'s they do, it computes their $\bx_s$'s and removes them from $\bTheta$. The algorithm iterates by randomly selecting another $\btheta$ until $\bTheta$ gets empty. We next comment on four implementation details: 

\emph{d1) Memory efficiency:} Each region is visited once and its parameters $(\bM,\br,\mcS_p\cap \mcS_d)$ are computed \emph{on the fly}, i.e., they are discarded once we have solved the OPF for each $\btheta$ of this region, thus yielding MP-OPF dramatic memory savings.

\emph{d2) Random sampling:} Ideally, one would like to visit regions in order of decreasing cardinality.
In this way, the dataset $\bTheta$ shrinks rapidly at Step 10, and thus, the number of checks at Step 8 drops significantly. Nonetheless, the cardinalities of regions are not known \emph{a priori}. In an attempt to visit regions with higher cardinality early on, we sample $\btheta$'s from $\bTheta$ randomly rather than following their order of appearance.

\emph{d3) Degeneracy:} MP-OPF introduces a region for $\btheta_o$ only if Step 4 checks affirmatively. Matrix $\bK$ can have rank deficiency due to linearly dependent constraints and/or primal degeneracy (more active constraints than variables). Although such cases could be handled, they would increase computational complexity. For this reason, we decided to only store their minimizers and not explore the related region. Our numerical tests demonstrate that such cases are relatively few.

\emph{d4) Identifying active constraints:} Step 3 identifies the active constraints of \eqref{eq:qp:ineq}, for which in theory the corresponding entries of vector $\be_o:=\bA\bx_o-\bE\btheta_o-\bb$ should be exactly zero. Interior point-based QP solvers though neither bring these entries to zero, nor they explicitly pinpoint the active constraints. To this end, we set solver's accuracy to $10^{-10}$ and deemed a constraint as active only when $\|\be_o\|_\infty\leq 10^{-5}$. For only $0.07\%$ of $\btheta$'s in $\bTheta$, our inferred partitioning of constraints into (in)active was incorrect. Such cases can be safely detected \emph{a posteriori}, since they did not satisfy $\btheta_o\in \mcS_p\cap\mcS_d$. For these few instances failing the latter sanity check, we only stored their minimizers and did not explore their regions. Finally, the check at Step 8 was deemed positive if all inequalities were satisfied with a precision of less than $10^{-4}$.

\begin{remark}\label{re:topology}
The proposed method assumes that the feeder topology remains fixed. Nevertheless, utilities oftentimes do reconfigure grids for the purposes of load balancing and loss minimization. With manual switches, a feeder may be reconfigured by line crew on a seasonal basis. In such cases, problem~\eqref{eq:qp} has to be solved for several fixed topologies. Topological changes alter matrices $(\bA,\bB,\bE,\bF,\bH)$ of \eqref{eq:qp} rather than its parameter vector $\btheta$. Since these matrices appear in products with the variable $\bx$, the KKT conditions of \eqref{eq:cond} become nonlinear in the problem parameters, and hence, the standard MPQP theory is not applicable. MP-OPF can still be adopted per topology to expedite PHCA. For remotely-controlled switches, their on/off status can be treated as variables of a possibly multistage OPF, where topology is optimized hourly and inverters are re-dispatched every 5~min or so. In this case, the inverter dispatch task would be a mixed-integer QP (MIQP). Although MPP-based approaches have been explored for MILP/MIQP~\cite{DBEP02}, \cite{DP09}, several of the neat properties of MPQPs are lost. Handling MIQP formulations for inverter dispatch that jointly optimize topologies goes beyond the scope of this work.
\end{remark}
\color{black}

\section{Numerical Tests}\label{sec:tests}
MP-OPF was tested on the IEEE 123-bus and a 1,160-bus feeder converted to single phase~\cite{jalali19}. Our PHCA was conducted using real-world hourly active load data extracted from the Iowa State University dataset~\cite{Iowa}. For solar generation data, we averaged the $15$~minute-based data provided by the Pecan Street dataset~\cite{pecandata} for 2017 to obtain their mean hourly values. The year-long load and solar sequences were randomly assigned to buses and scaled so their peak values matched the nominal load values. Lacking records for load power factors, reactive loads were synthesized by independently drawing lagging power factors at random within $[0.90,0.95]$ for each bus and keeping them fixed throughout the year. 

\cmb{We solved \eqref{eq:opfr} for $\beta=0.2$ for the 123-bus system and for $\beta=0$ for the 1,160-bus system.} To improve the numerical conditioning of \eqref{eq:qp}, we scaled $(\bH,\bC,\bd)$ of \eqref{eq:qp:cost} by $\|\bH\|_2$; $(\bA,\bE,\bb)$ of \eqref{eq:qp:ineq} by $\|\bA\|_2$; and $(\bB,\bF,\bf)$ of \eqref{eq:qp:eq} by $\|\bB\|_2$. To select $\eta$, we solved $1,000$ random feasible instances of \eqref{eq:opf} and recorded their Lagrange multipliers. We set $\eta=1$ to satisfy Proposition~\ref{prop:exact penalty} for the tested instances, and $\nu=20$, the largest eigenvalue of $\bH$. Problem~\eqref{eq:qp} was solved on a computer with Intel Core7 @ 3.4 GHz (16 GB RAM) with MATLAB using YALMIP~\cite{YALMIP}, SeDuMi~\cite{sedumi}, and ECOS~\cite{ECOS}.

For the IEEE 123-bus system, our PHCA involves $3$ analysis parameters: \emph{i) Injection scaling:} (re)active injections were scaled by $\{1,2,3\}$ times their nominal values to also capture overloaded conditions; \emph{ii) Inverter oversizing:} To study the necessity of reactive power compensation at peak solar irradiance, inverters were deployed  with $10\%$ oversizing $(\bar{s}_n^g=1.1\bar{p}_n^g)$ and without oversizing $(\bar{s}_n^g=\bar{p}_n^g)$ for each bus $n$~\cite{Turitsyn11}; and \emph{iii) DER penetration levels} varying from $10\%$ to $100\%$ at increments of $10\%$. These $3$ analysis parameters were used to construct $3\times 2 \times 10=60$ different values of $\mcA$. The uncertain parameter set $\mcU$ consisted of hourly load/solar tuples across all buses for a period of $1$ year, that is $8,640$ members for $\mcU$. Combining $\mcA$ with $\mcU$ resulted in a $\bTheta$ with $S=|\mcA|\times |\mcU|=518,400$ instances.

To analyze the scalability of the proposed PHCA algorithm, our MPP-aided PHCA algorithm was further tested on a 1,160-bus system. This system was obtained from the IEEE 8,500-bus feeder upon maintaining the buses located under the voltage regulator of bus 338. This system is approximately ten times larger than the IEEE 123-bus feeder. For zero-injection buses, we added independent identically distributed Gaussian random power ratings of mean $p_{min}$ and variance $0.01$ with $p_{min}$ selected as the smallest power rating among non-zero injection buses. The rest of the data generation process was similar to that for the IEEE $123$-bus feeder. 

For the tests on the 1,160-bus feeder, DER levels varied from $10\%$ to $100\%$ with increments of $10\%$; inverters were oversized by $10\%$; and power injections were not scaled. Hence overall, the analysis parameters entailed $10$ different options. Combining $\mcA$ and $\mcU$ resulted in a $\bTheta$ with $86,400$ instances. The analysis parameters entailed $10$ different options. Combining $\mcA$ and $\mcU$ resulted in a $\bTheta$ of 86,400 instances. We did not process the entire IEEE 8,500-bus feeder, because the aforementioned dataset $\bTheta$ would not fit in the RAM (random access memory) of the personal computer used for these tests. To circumvent that difficulty, one should potentially modify our MATLAB code and access one record of $\bTheta$ at a time with appropriate read/write commands. MP-OPF is amenable to such implementation for the scanning process at Step 7. Such serial data access would slow down the overall solution time, but that would is unavoidable for any other solution methodology anyway. We did not pursue this access and considered the 1,160-bus reduced version instead.

\begin{table}

	\renewcommand{\arraystretch}{1.2}
	\caption{Mean and Standard Deviation (std) of QP Solver Times for OPF [sec]}
	\label{tbl:time} \centering
	\begin{tabular}{|l|r|r|r|r|r|r|}
		\hline\hline
	    \multirow{2}{*}{\textbf{Feeder}} &   \multicolumn{2}{c|}{\textbf{YALMIP}} & \multicolumn{2}{c|}{\textbf{SeDuMi}}  & \multicolumn{2}{c|}{\textbf{ECOS}}\\
	   \cline{2-7}
	    &  mean & std & mean & std  & mean & std\\
		\hline\hline 123-bus&$0.100$ & $0.005$ & $0.244$ & $0.019$ & $0.126$ & $0.012$\\
		\hline	1,160-bus&$1.092$ & $0.104$ &$16.680$ & $1.150$ & $25.890$ & $1.230$\\
		\hline\hline
	\end{tabular}
	\color{black}
\end{table}

\begin{figure}[t]
		\centering
	\subfigure{\includegraphics[scale=0.6]{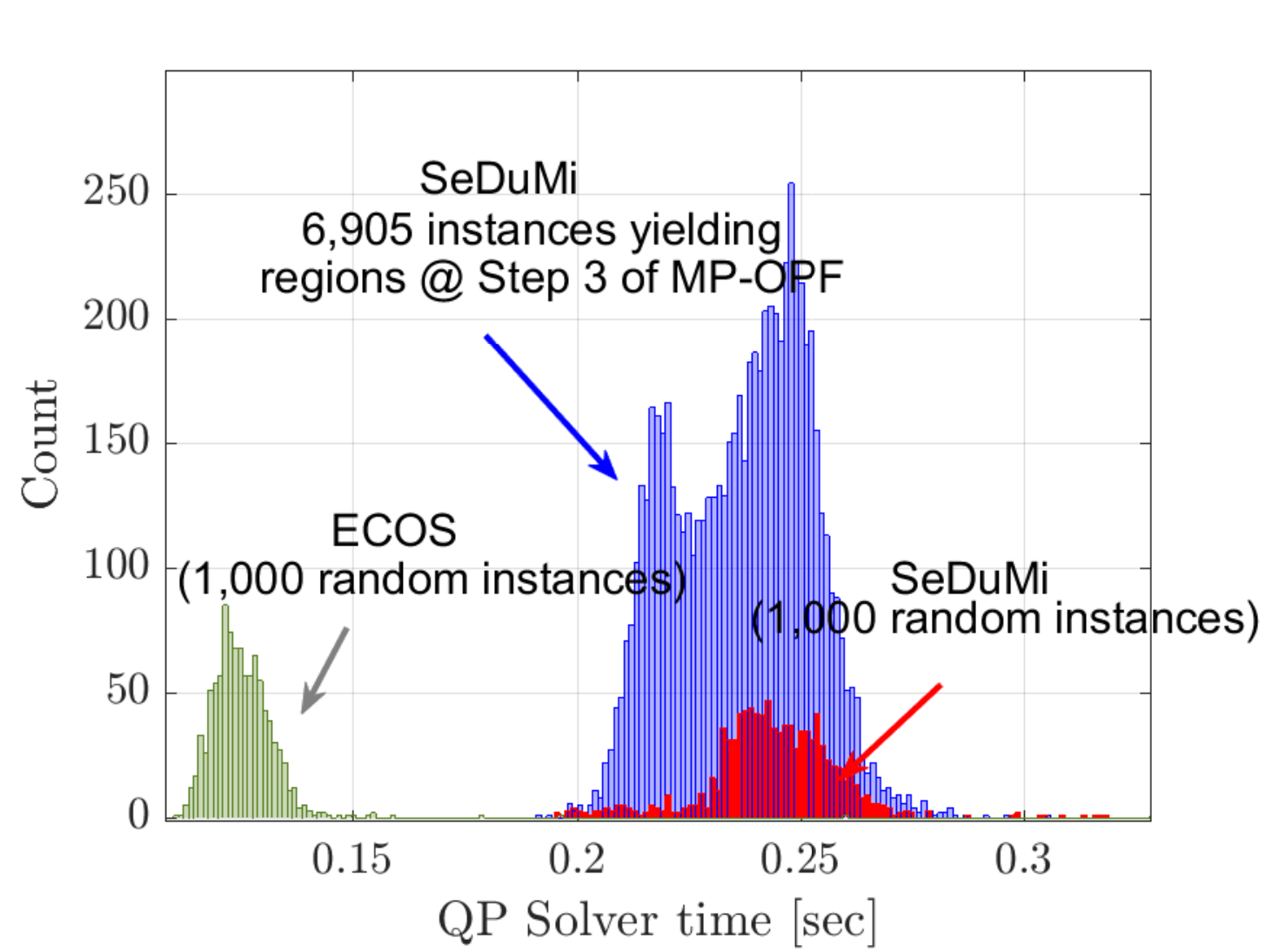}}
		\caption{The histograms of QP running times for SeDuMi and ECOS over 1,000 randomly drawn $\btheta$'s for the IEEE 123-bus system. The histogram of the QP running times solved by SeDuMi during Step 3 of MP-OPF are also shown for the 6,905 instances that yielded a region.}
	\label{fig:hist_time}
\end{figure}

MP-OPF was compared against solving a QP for each one of the $S$ $\btheta$'s. Since the latter is computationally taxing, we solved \eqref{eq:qp} only for $1,000~\btheta$'s randomly drawn from $\bTheta$ with SeDuMi and ECOS. Then, we used the average solution times to extrapolate over all OPF instances, and thus estimate the time those solvers would require to solve PHCA. Figure~\ref{fig:hist_time} depicts the histograms of the recorded times over these 1,000 randomly sampled OPF instances. It also shows the histogram for SeDuMi running times experienced while solving the 6,905 OPF instances at Step 3 of MP-OPF that generated a critical region. Table~\ref{tbl:time} shows the mean and standard deviation of these solution times. Times are reported individually for YALMIP translation and the respective solver (SeDuMi or ECOS) alone. 

From Figure~\ref{fig:hist_time} and Table~\ref{tbl:time}, ECOS is faster than SeDuMi for the 123-bus system, but slower than SeDuMi for the 1,160-bus system. Because of that, when PHCA was to be solved without the aid of MPP, we used ECOS to solve the QPs for the 123-bus system and SeDuMi to solve the QPs for the 1,160-system system. For either systems however, Step~3 of our MP-OPF ran SeDuMi (despite not always being the faster solver). This is because we wanted higher accuracy in determining whether a constraint is active or not, so regions can be identified and represented reliably. It is for this reason that Figure~\ref{fig:hist_time} shows only the running times of SeDuMi over Step~3 of MP-OPF, and not ECOS.

It is also worth mentioning from Table~\ref{tbl:time} that the standard deviations are much smaller than the mean values for running times, implying that the QP solution time is relatively insensitive to the particular value of $\btheta$. Because of that, rather than solving each QP of the corresponding dataset, we can safely estimate the total running time for PHCA without MP-OPF by multiplying the mean values of Table~\ref{tbl:time} by the total number of instances. Based on this estimate, the conventional PHCA would take roughly $49.6$~hours using SeDuMi for the $123$-bus system and $408$ hours for the $1,160$-bus system. The same analysis with the ECOS solver would have taken $32.5$ and $650.8$ hours, respectively. Note that all reported times are \emph{wall-clock} times.

\begin{figure}[t]
	\centering
	\subfigure{\includegraphics[width=0.38\textwidth]{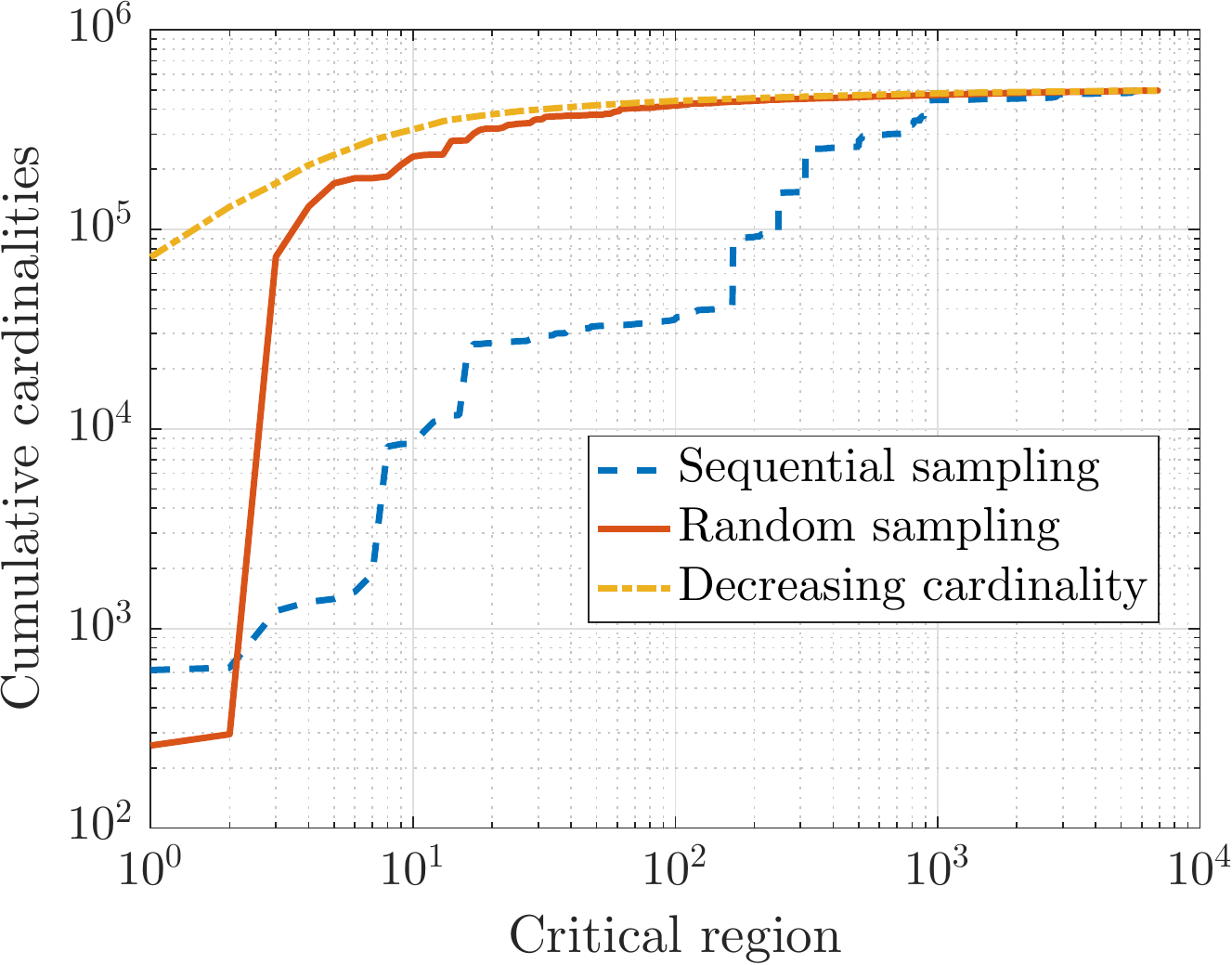}}
	\caption{On the horizontal axis, regions are indexed based on: \emph{i)} the order they appear in $\bTheta$; \emph{ii)} the order they were visited by randomly drawing $\btheta$'s; and \emph{iii)} based on decreasing cardinality (most populous are visited first). The vertical axis shows the cumulative sum of region cardinalities. For example, from the top curve we see that the largest region includes roughly 70,000 instances, the two largest regions include 110,000 instances, and so on. Ideally, MP-OPF would perform best if regions were visited in decreasing cardinality. Tnis is not practical since region cardinalities are not known \emph{a priori}. Nonetheless, by randomly sampling instances from the dataset rather than accessing them in a sequential (chronological) order, we can get closer to the ideal case and reduce computational time from 4.6 to 3.2 hours.}
	\label{fig:rhist}
\end{figure}

\begin{figure}[t]
	\vspace*{1em}
	\centering
	\subfigure{\includegraphics[width=0.39\textwidth]{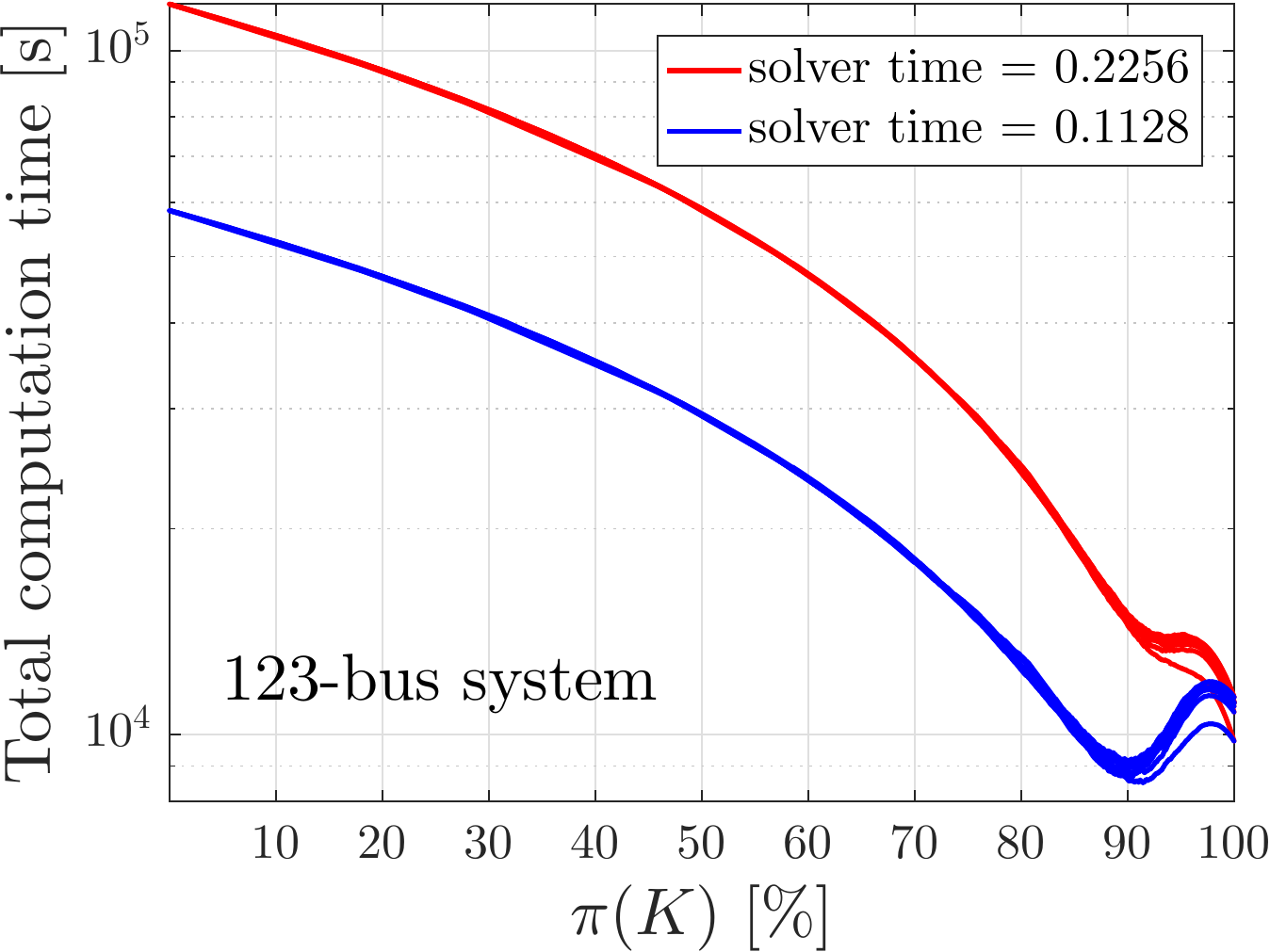}}
	\subfigure{\includegraphics[width=0.39\textwidth]{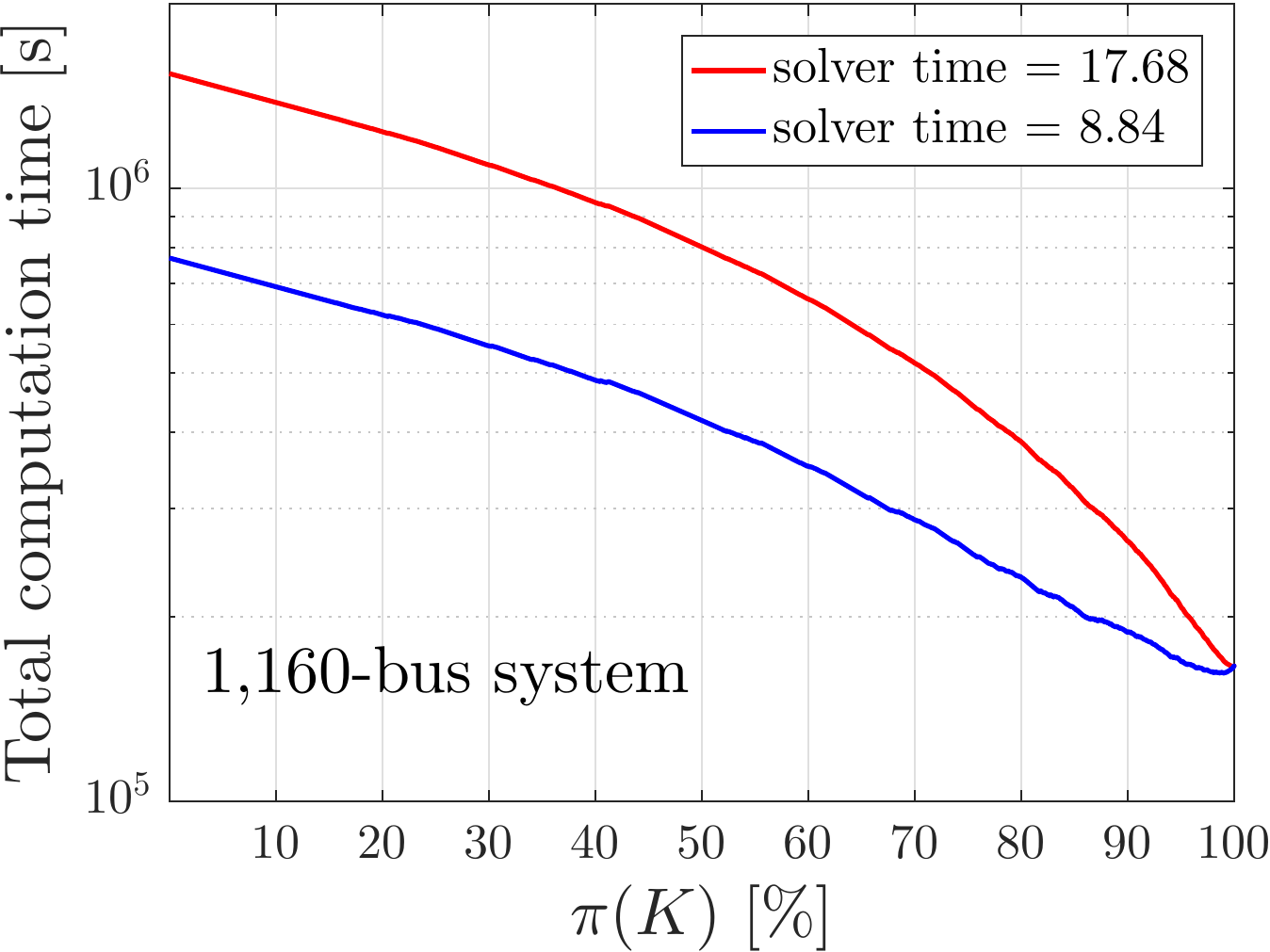}}
	\caption{Monte-Carlo tests to explore the potential benefit of early stopping MP-OPF depending on the solution time of the QP solver. \emph{Top panel:} For the 123-bus system and ECOS's solution time of $0.2256$~sec per QP, apparently there is no benefit. For a hypothetical faster solver with half the solution time of ECOS, there is a benefit of stopping after $91\%$ of $\btheta$'s are handled by MP-OPF. \emph{Bottom panel:} For the 1,160-bus system and SeDuMi's solution time of $17.68$~sec per QP, there is no benefit in early stopping even using a hypothetical solver with half SeDuMi's time.}
	\label{fig:10tests}
\end{figure}

When the same task was solved with MP-OPF under SeDuMi, we visited $6,905$ regions for the $123$-bus system and it took only $3.2$~hours. For the $1,160$-bus system, we visited $2,111$ regions and it took $42.6$ hours. This is an improvement by an order of magnitude over the conventional approaches. For the $123$-bus system, when the random sampling of $\btheta$'s in Step 2 was replaced by sequential sampling (visiting $\btheta$'s in the chronological order they appear in $\bTheta$), MP-OPF took $4.6$~hours. This advantage of random versus sequential sampling is explained in Figure~\ref{fig:rhist}. Its curves depict that randomly sampling $\btheta$'s entails visiting critical regions in almost decreasing order. This is reasonable since it is more likely to draw $\btheta$'s from more populous regions. The top curve also shows that more than $95\%$ ($90\%$) of $\btheta$'s belong to only $500$ ($80$) critical regions, which is a testament to the timing advantage of MP-OPF.

Spurred by the online identification of active sets based on the rate of discovery introduced in~\cite{MRN19} and the observation that relatively few regions contain the majority of OPF instances (Figure~\ref{fig:rhist}), we next examined if there is any benefit in not exploring every single critical region. We tried terminating MP-OPF for both system prematurely after exploring only $K$ critical regions. The unclaimed $\btheta$'s were dealt with by solving the individual OPFs one by one. Number $K$ was varied from $1$ to the total number of regions. For each $K$, we recorded the time $t_m(K)$ spent on MP-OPF and the percentage $\pi(K)$ of $\btheta$'s handled by the algorithm. Since the solution times by both solvers were shown earlier to be rather invariant to the particular $\btheta$, we assumed an average solution time of $0.2256$~sec for the $123$-bus system and $17.68$~sec for the $1,160$-bus system per instance and estimated the time needed for handling the unclaimed $\btheta$'s accordingly as $t_s(K)$. Apparently $t_m(K)$ is increasing and $t_s(K)$ is decreasing with $K$. For the $123$-bus system, this test was repeated for 10 Monte-Carlo trials over randomly shuffled datasets. Figure~\ref{fig:10tests} (top panel) depicts the total running time over $\pi(K)$. The figure suggests that there is no benefit in not exploring all regions. The explanation behind this is that Steps 4-10 of MP-OPF are computationally inexpensive compared to the QP solver, especially when $\bTheta$ becomes smaller. To explore the effect of the QP solver's speed on the previous discussion, we repeated the previous test presuming a hypothetical QP solver with half the solution time ($0.1128$~seconds for the $123$-bus and $8.84$~sec for the $1,160$-bus systems) per QP instance. The timing results of Fig.~\ref{fig:10tests} (bottom panel) demonstrate that for the $123$-bus system and with a faster solver, it is computationally advantageous to stop MP-OPF after $91\%$ of $\btheta$'s have been covered. For the $1,160$ system, there is no such benefit even with a faster QP solver.
\color{black}

The next test intends to show that the amortized advantage of MP-OPF increases as the OPF dataset $\bTheta$ grows bigger. To this end, we compared the timing result of the previous test for the $123$-bus system with the PHCA setup where injection scaling took only one of the values $\{1,2,3\}$. The new set $\mcA$ involved $1\times 2 \times 10=20$ options. When those options are combined with $\mcU$ give a smaller dataset of $172,800$ instances. The times for the individual injection scalings of $\{1,2,3\}$ was $\{0.23,1.04,3.06\}$~hours respectively, yielding the combined duration of $4.3$~hours. This indicates that MP-OPF becomes more efficient as the cardinality of regions grows and so region descriptions are reused.

\begin{figure*}[t]
	\centering	
	\subfigure{\includegraphics[width=0.32\textwidth]{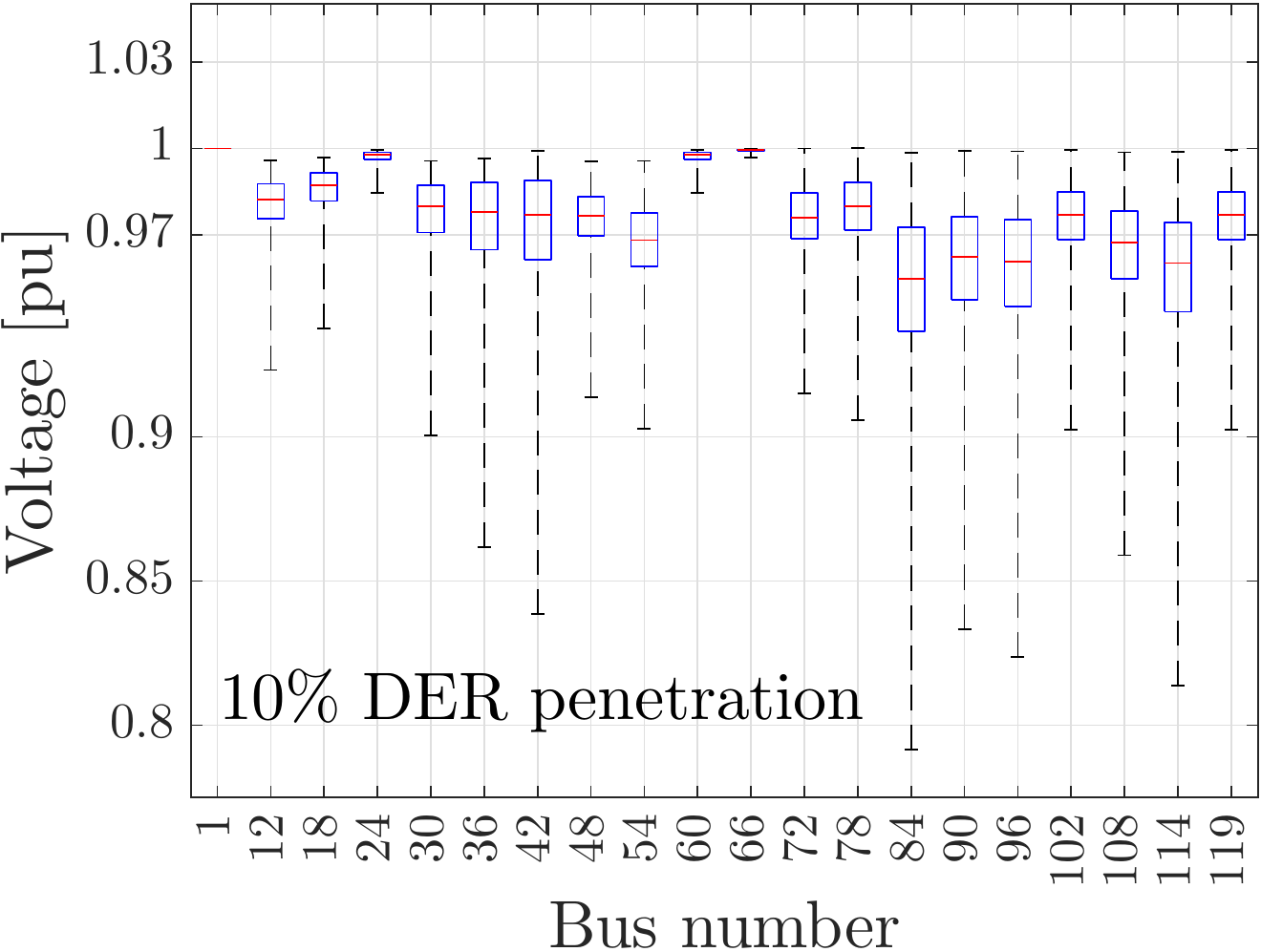}}
	\hspace*{+1em}
	\subfigure{\includegraphics[width=0.305\textwidth]{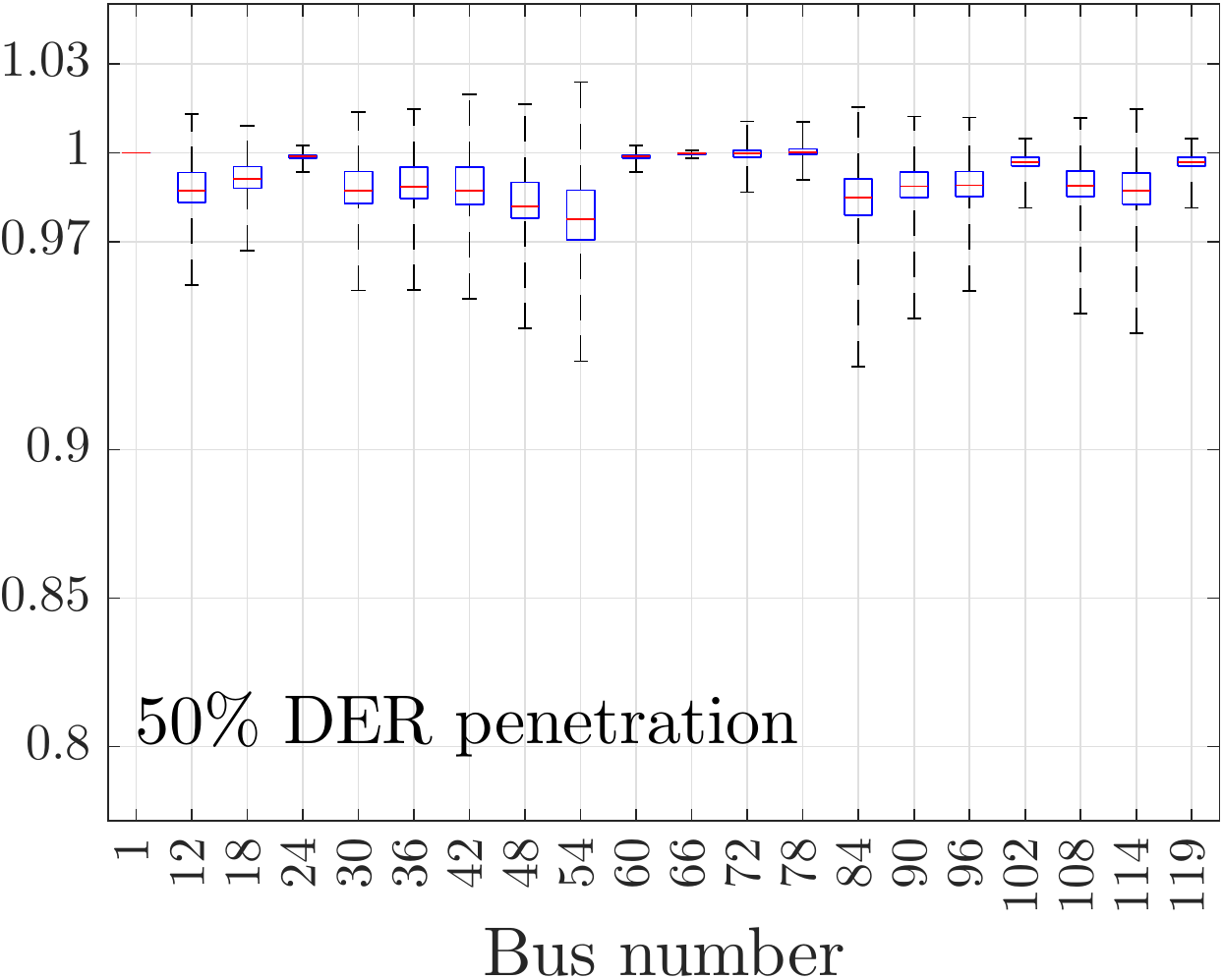}}
	\hspace*{+1em}
	\subfigure{\includegraphics[width=0.305\textwidth]{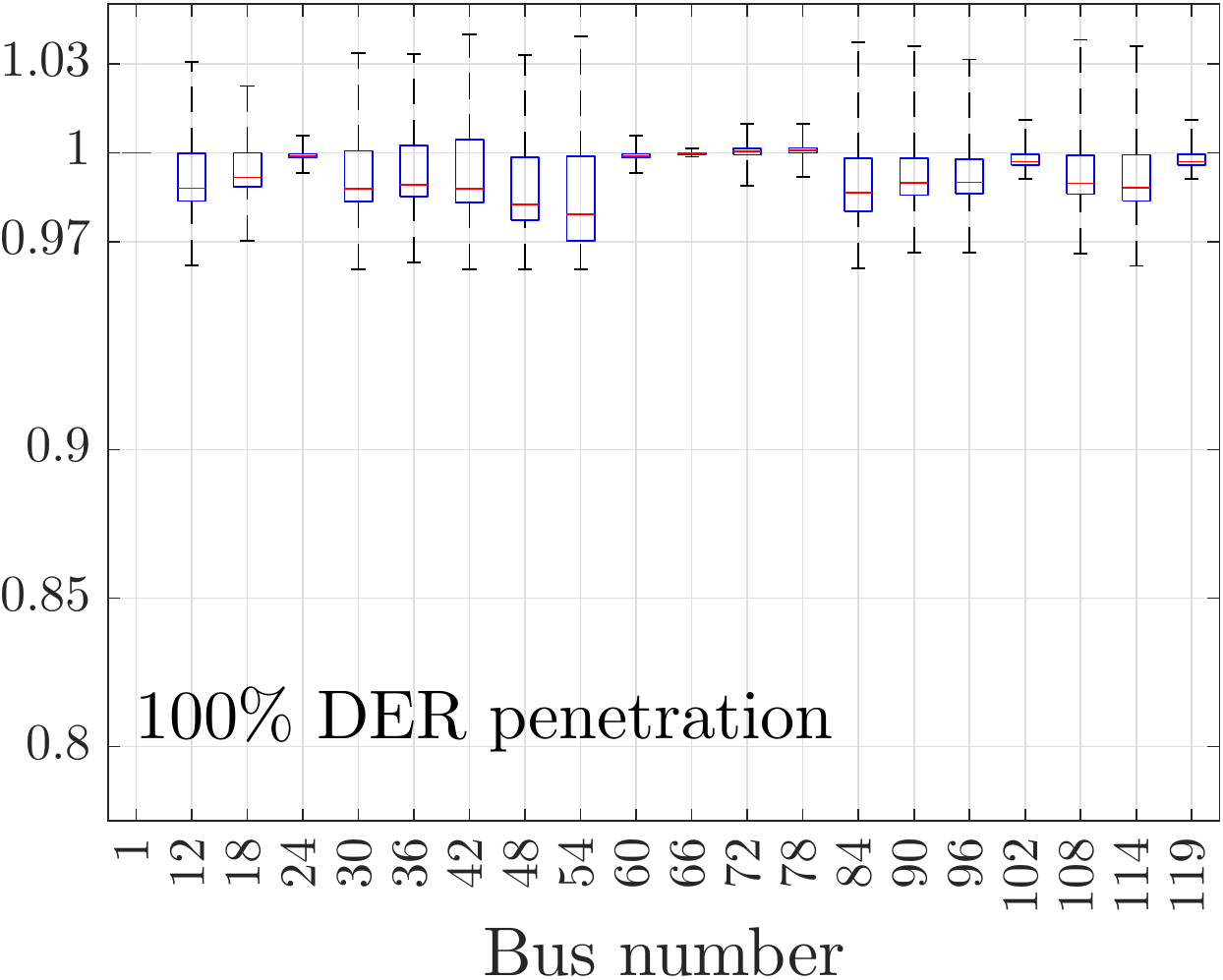}}
	\caption{Box plots for a selection of bus voltages for the injection scaling of $3$. Bus 1 is the substation, while ${24,60,66}$ host LDC regulators.}
	\label{fig:boxplots}
\end{figure*}

\begin{figure*}[t]
	\centering	
	\subfigure{\includegraphics[width=0.32\textwidth]{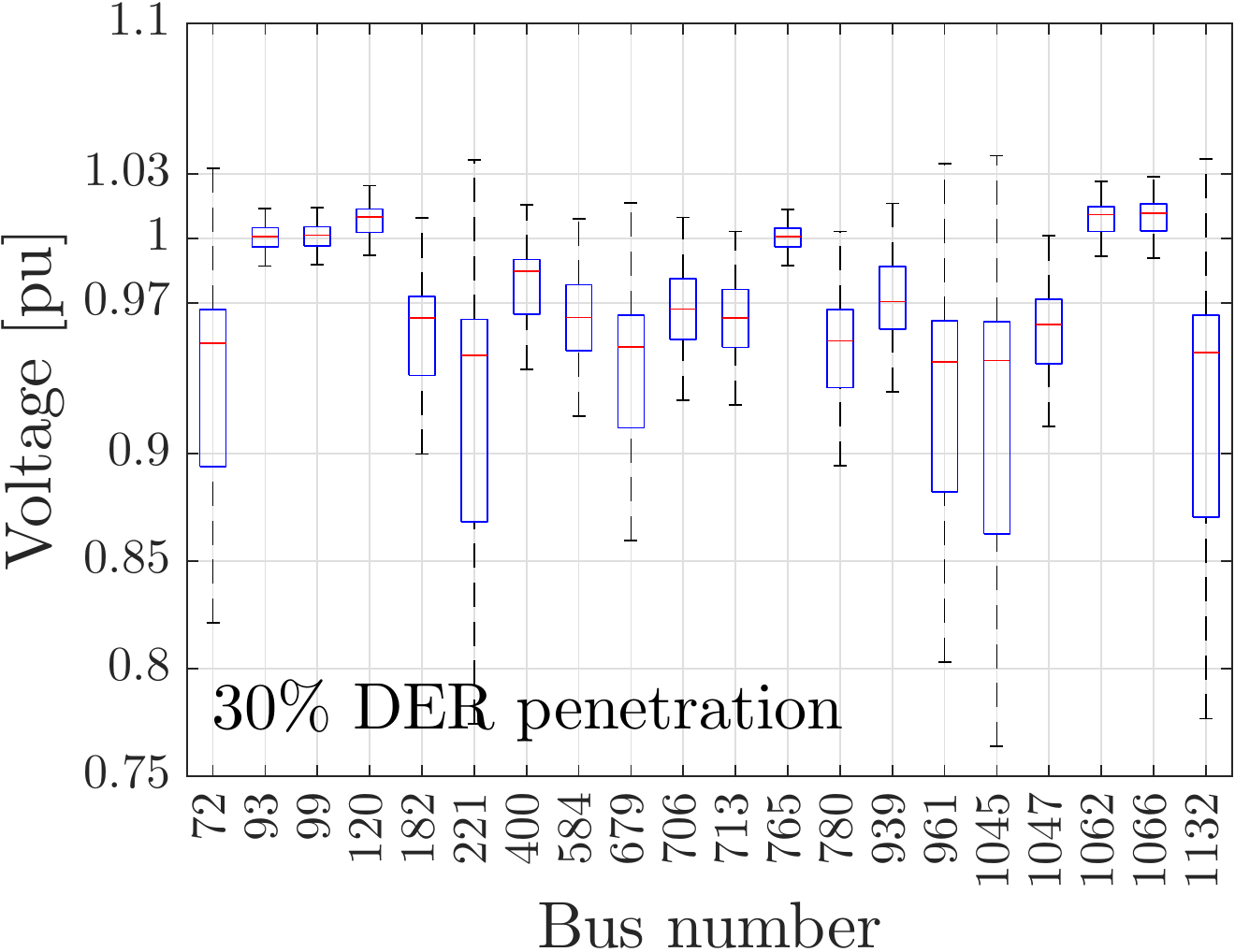}}
	\hspace*{+1em}
	\subfigure{\includegraphics[width=0.305\textwidth]{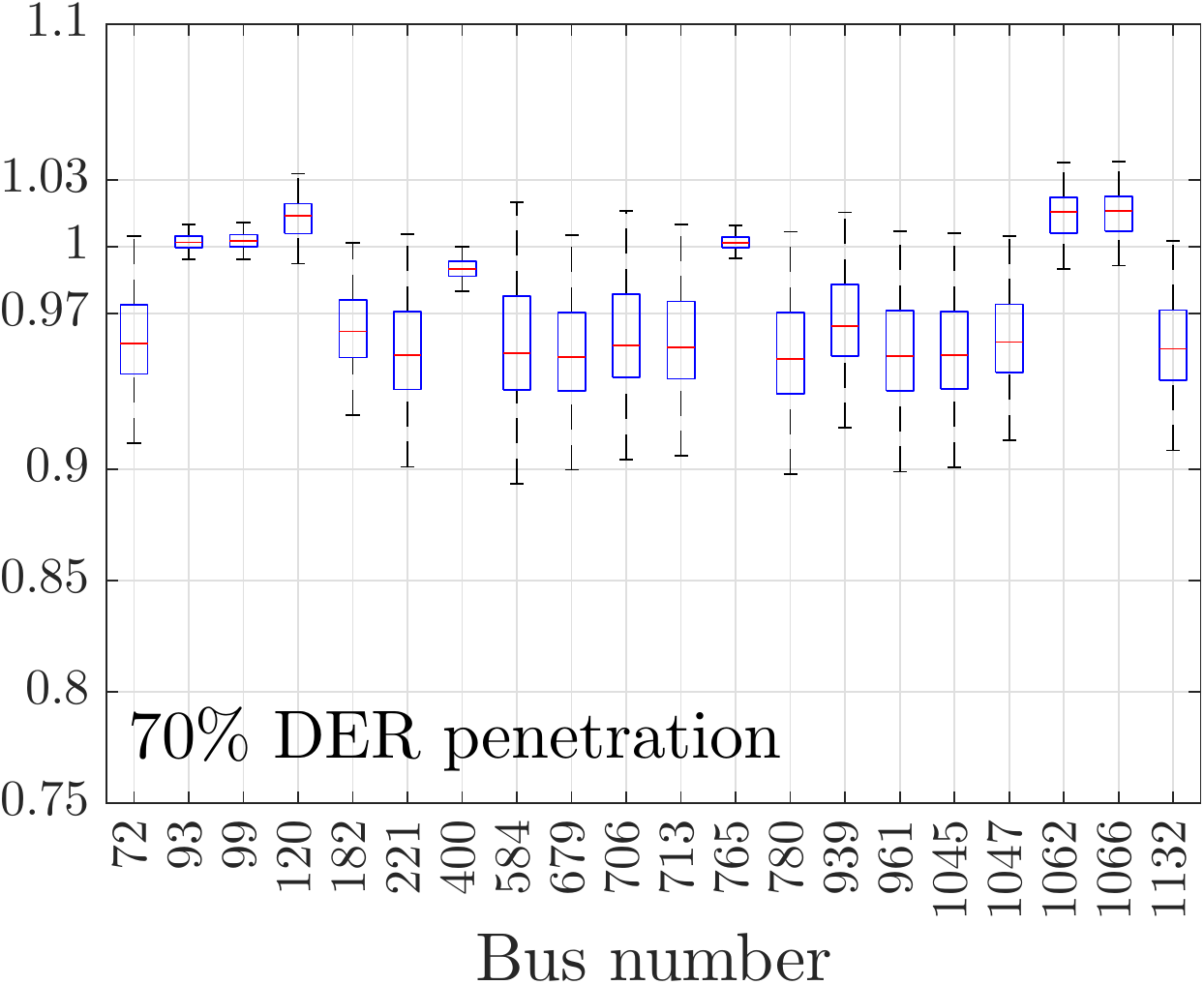}}
	\hspace*{+1em}
	\subfigure{\includegraphics[width=0.305\textwidth]{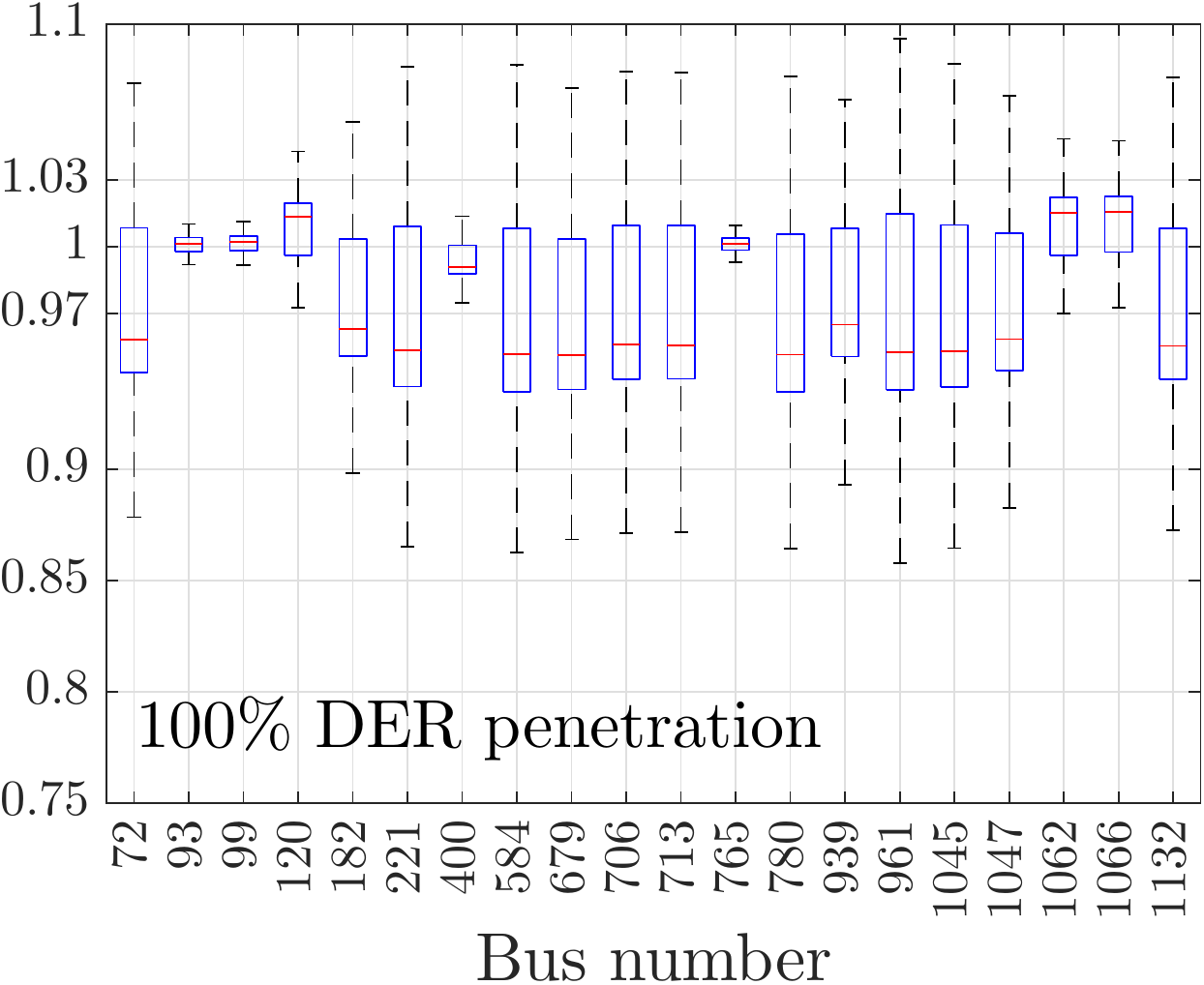}}
	\caption{Box plots for a selection of bus voltages for the 1,160-bus system for three solar penetration levels.}
	\label{fig:boxplots8500}
\end{figure*}

\begin{figure}[t]
	\centering
	\vspace*{-0.3em}
	\subfigure{\includegraphics[width=0.36\textwidth]{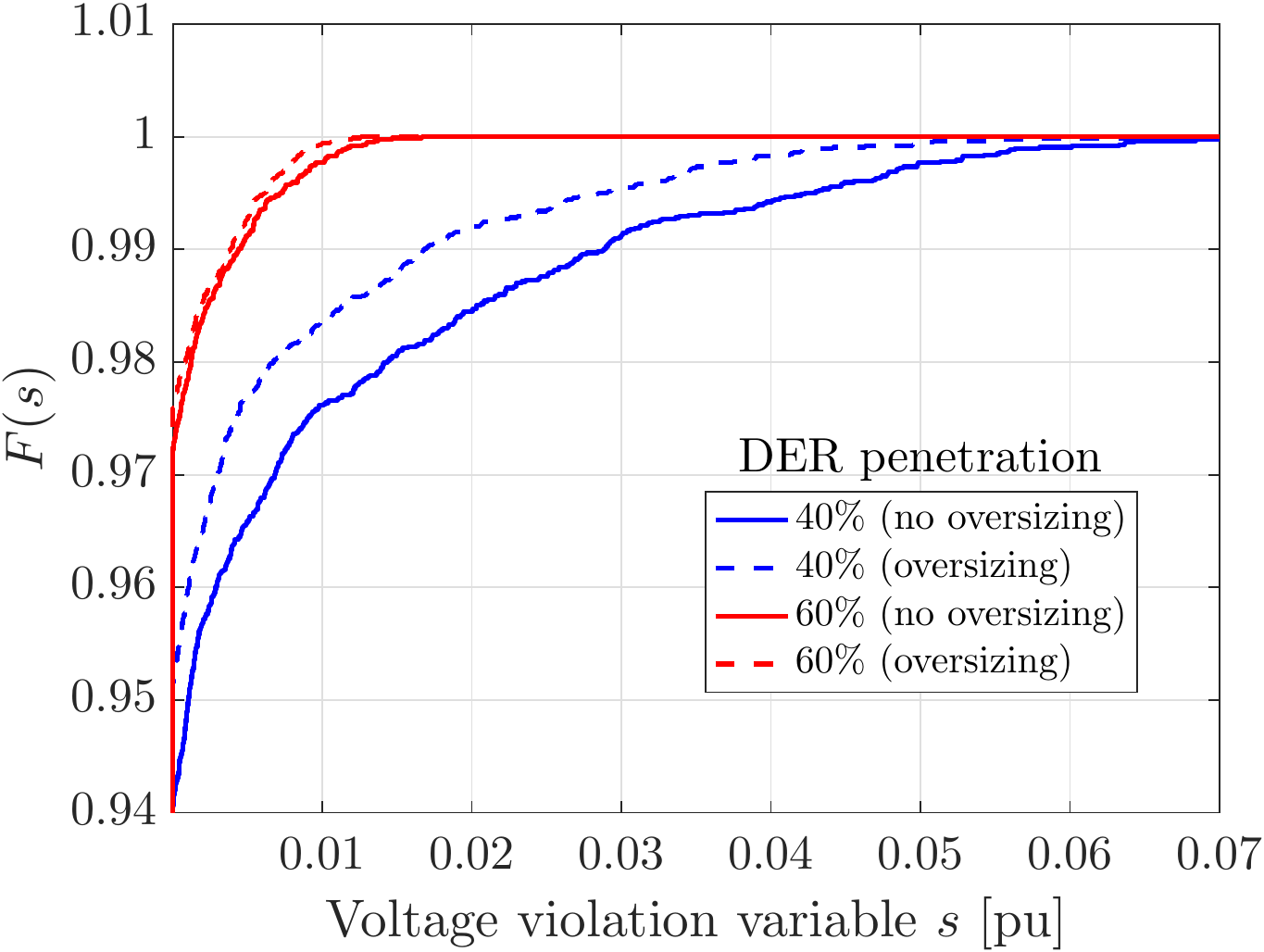}}
    \subfigure{\includegraphics[width=0.35\textwidth]{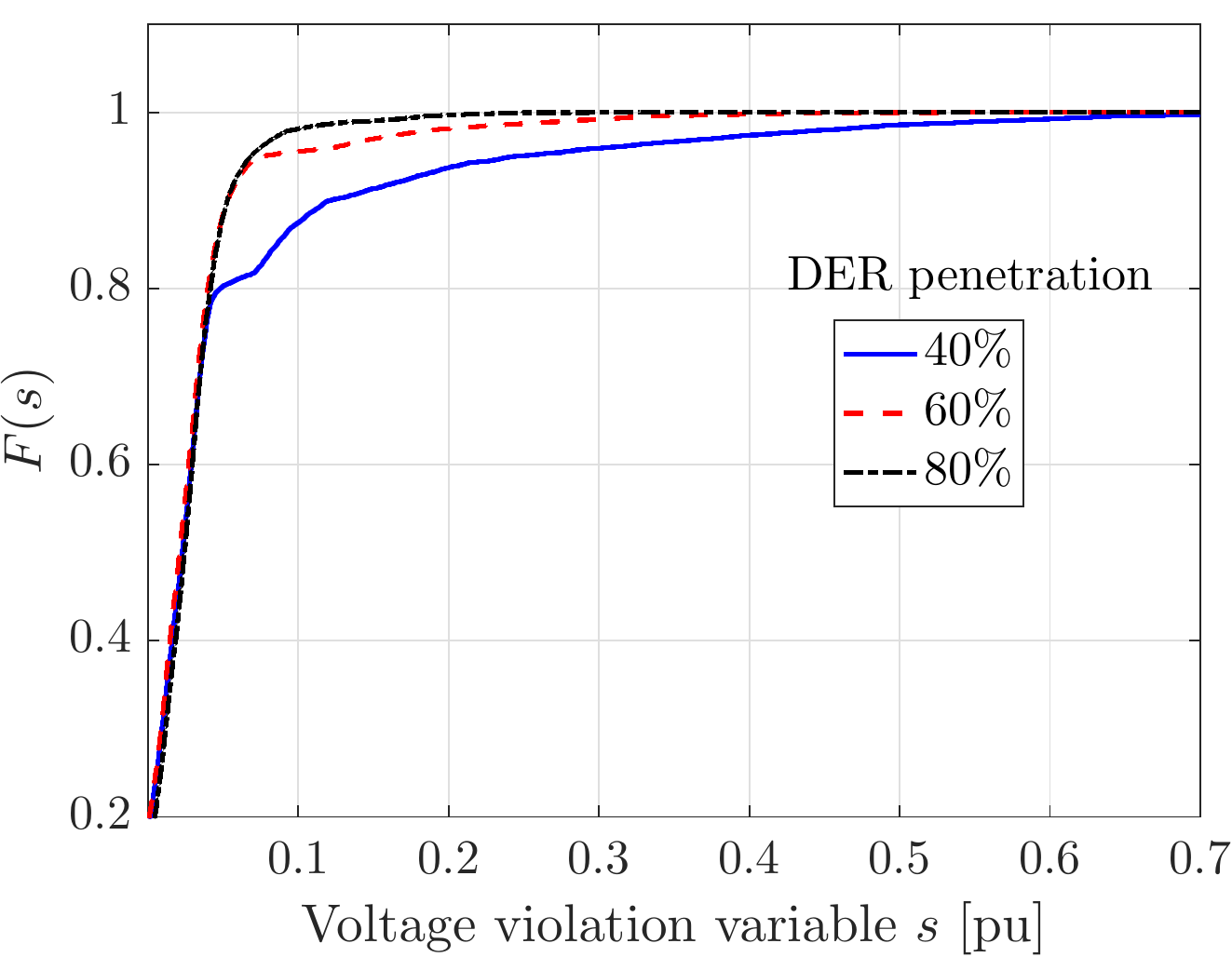}}
	\caption{\cmb{The cdf of $s$ for two DER penetration levels during 10am-3pm for IEEE 123-bus system (top) and for all hours for the 1,160-bus network (bottom).}}
	\label{fig:CDF-scales}
\end{figure}

Let us now showcase a few of the grid statistics one can readily compute with PHCA through MP-OPF. Figure~\ref{fig:boxplots} (top) presents the distribution of voltage violations for the $123$-bus system when injections are scaled up by a factor of $3$. Under-voltages occur at the $10\%$ DER penetration level. These problems are alleviated by reactive power control from DERs at $50\%$ and $100\%$ penetrations. On the other hand, at the $100\%$ penetration, over-voltages are likely to occur. For the 1,160-bus system, Figure~\ref{fig:boxplots} (bottom) plots voltage statistics for three solar penetration levels. At $30\%$ penetration level, the feeder mainly experiences under-voltage problems. By increasing the penetration level to $70\%$, under-voltage problems are alleviated to some degree. When solar penetration level is increased to $100\%$, more over-voltage problems occur.

We also studied the overall voltage violation statistics. Figure~\ref{fig:CDF-scales} (top) depicts the cdf for the constraint violation variable $s$ under injection scaling by $3$. The solid blue (resp. red) curve shows that at $40\%$ (resp. $60\%$) DER penetration, voltages remain within limits for $94\%$ (resp. $97\%$) of the time and the maximum constraint violation is upper bounded by $0.07$~pu ($0.02$~pu). This improvement with increasing solar penetration is attributed to DER control. The curves also corroborate the benefit of inverter oversizing for enhanced voltage regulation. Figure~\ref{fig:CDF-scales} (bottom) depicts the cdf of the constraint violation variable $s$ for the $1,160$-bus system for three penetration levels. This network suffers from higher constraint violations with respect to the $123$-bus system. However, with higher solar penetration, these violations happen with less frequency and severity. 
\color{black}

\section{Conclusions}
Viewed from the vantage point of MPP, the formidable task of PHCA has been expedited by a factor of 10. Rather than solving a large number of OPF instances, MPP allows the utility to actually solve only a small fraction of OPFs. The computational advantage is that the minimizers for the remaining OPFs are readily computed in closed form thanks to the parametric dependence of the OPF on analysis and uncertain parameters. The proposed PHCA relies on historical scenarios to swiftly infer sample statistics for all grid quantities of interest at once, and without knowing the probability distribution for load/solar. To the best of our knowledge, this is the first time MPP is used in distribution grid operations. Albeit framed within the context of PHCA, MPP can be beneficial in speeding up DER inverter control. Yet in this case, MP-OPF has to be adapted to operate in an online instead of a batch setup. Extending MP-OPF to multiphase grids is a relatively straightforward extension, whereas using it towards AC-OPF problems is a more challenging research direction. 

\appendices

\section{Proof of Proposition~\ref{pro:losses}}\label{sec:appendixA}
Let vector $\bphi$ collect the angles of voltage phasors $V_n=v_ne^{j\phi_n}$ across all buses of the feeder excluding the substation bus. Define also the vectors of power injections $\bs:=[\bp^\top~\bq^\top]^\top$ and voltages in polar coordinates $\bu:=[\bv^\top~\bphi^\top]^\top$. Consider the distribution line running from bus $m$ to bus $n$ as introduced before \eqref{eq:Vmag}. Because this line feeds bus $n$, we can simply index this line by $n$; recall that each bus is fed by a single line due to the radial structure of the feeder. By definition, the ohmic losses on this line are
\[L_n(\bu)=g_{mn}\left(v_m^2+v_n^2-2v_mv_n\cos(\psi_n)\right)\]
where $g_{mn}=r_{mn}/(r_{mn}^2+x_{mn}^2)$ is the line conductance and $\psi_n:=\phi_m-\phi_n$. At the flat voltage profile $\bbu:=[\bone^\top~\bzero^\top]^\top$, we apparently get $L_n(\bbu)=0$.

We next compute the first- and second-order partial derivatives of $L_n$ with respect to $\bu$. Let $\bu_{n}:=[v_m~v_n~\phi_m~\phi_n]^\top$ be the subvector of $\bu$ collecting the voltages appearing at the ends of line $n$. It is easy to compute the gradient vector
\[\nabla_{\bu_n}L_n(\bu)=2g_n\left[\begin{array}{l}
v_m-v_n\cos\psi_n\\
v_n-v_m\cos\psi_n \\
+v_m v_n\sin\psi_n\\
-v_mv_n\sin\psi_n
\end{array}
\right].\]
The remaining entries of $\nabla_{\bu_n}L_n(\bu)$ are obviously zero. Evaluating this gradient vector at $\bbu$, we get $\nabla_{\bu}L_n(\bbu)=\bzero$. We can similarly compute the second-order partial derivatives of $L_n(\bu)$. Evaluating the related Hessian matrix at $\bbu$, we get
\begin{equation}\label{eq:hessian}
\nabla_{\bu_n\bu_n}^2L_n(\bbu)=2g_n\left[\begin{array}{cccc}
+1 &   -1 &   0 &   0\\
-1 &   +1   &   0   &   0\\
0 &   0   &   +1   &   -1\\
0 &  0   &   -1   &   +1
\end{array}
\right].
\end{equation}

The remaining entries of the complete Hessian matrix $\nabla_{\bu\bu}^2L_n(\bbu)$ are apparently zero. To simplify notation, let us henceforth denote $\nabla_{\bu\bu}^2L_n(\bbu)$ as $\tbN_n$. If $\bgamma_n$ is the $n$-th row of the reduced branch-bus incidence matrix $\bGamma$, then the Hessian matrix $\tbN_n$ can be compactly written as
\begin{equation}\label{eq:tbNn}
\tbN_n=\left[\begin{array}{cc}
\bN_n & \bzero\\
\bzero & \bN_n
\end{array}\right]\quad \text{where}\quad \bN_n:=2g_n\bgamma_n\bgamma_n^\top.
\end{equation}

We would like to express the losses on line $n$ as a function of $\bs$, i.e., $L_n(\bs)$, not as a function of $\bu$. Since line losses are a complicate function of power injections, we resort to a Taylor's series approximation around the flat voltage profile. Note that at $\bu=\bbu$, we also get zero power injections $\bbs=\bzero$ and zero line losses $L_n(\bbs)=0$. Aiming for a first-order approximation of $L_n(\bs)$, we can apply the chain rule to get
\[L_n^{(1)}(\bs)\simeq L_n(\bbs) + (\nabla_{\bu}L_n(\bbu))^\top\cdot \nabla_{\bs}\bbu\cdot (\bs-\bbs)\]
where $\nabla_{\bs}\bbu$ denotes the Jacobian matrix of $\bu$ with respect to $\bs$ evaluated at $\bbu$. Because $\nabla_{\bu}L_n(\bbu)=\bzero$ as argued earlier, we obtain $L_n^{(1)}(\bs)=0$ for all $\bs$. 

A second-order expansion of $L_n(\bs)$ around $\bbs=\bzero$ yields
\[L_n^{(2)}(\bs)\simeq \frac{1}{2}\bs^\top\bH_n\bs\]
where $\bH_n$ is the Hessian matrix $\nabla_{\bs\bs}^2L_n(\bs)$ evaluated at $\bs=\bbs$. The $(i,j)$-th entry of $\nabla_{\bs\bs}^2L_n(\bs)$ can be computed as
\begin{equation}\label{eq:Hij}
\frac{\partial^2 L_n}{\partial s_i\partial s_j}=\sum_{k,\ell}\frac{\partial^2 L_n}{\partial u_k\partial u_\ell}\cdot \frac{\partial u_\ell}{\partial s_i} \cdot\frac{\partial u_k}{\partial s_j}+
\sum_{k}\frac{\partial L_n}{\partial u_k}\cdot \frac{\partial^2 u_k}{\partial u_i \partial s_j}
\end{equation}
where we have used the three ensuing properties:
\begin{align*}
&\frac{\partial L_n}{\partial s_j}=\sum_{k} \frac{\partial L_n}{\partial u_k}\cdot \frac{\partial u_k}{\partial s_j}\\
&\frac{\partial}{\partial s_i}\left(\frac{\partial L_n}{\partial u_k}\right)=\sum_{k} \frac{\partial}{\partial s_i}\left(\frac{\partial L_n}{\partial u_k}\right)\cdot \frac{\partial u_k}{\partial s_j} + \sum_{k} \frac{\partial L_n}{\partial u_k}\cdot \frac{\partial^2 u_k}{\partial s_i \partial s_j}\\
& \frac{\partial}{\partial s_i}\left(\frac{\partial L_n}{\partial u_k}\right)=\sum_{\ell}  \frac{\partial^2 L_n}{\partial u_k \partial u_\ell}\cdot \frac{\partial u_\ell}{\partial s_i}.
\end{align*}

If we evaluate \eqref{eq:Hij} at $\bbs$, the second summation disappears because $\nabla_{\bu}L_n(\bbu)=\bzero$ as shown earlier. To evaluate the first summation, we need two components:
\renewcommand{\labelenumi}{\emph{c\arabic{enumi})}}
\begin{enumerate}
\item The second-order partial derivatives of $L_n$ with respect to voltages. These have been derived in \eqref{eq:hessian}.
\item The Jacobian matrix $\nabla_{\bs}\bu$ evaluated at $\bu=\bbu$. References \cite{BoDo15}, \cite{Deka1} have obtained a first-order Taylor's series expansion of $\bu$ with respect to $\bs$ at $\bbs$, which reads
\begin{equation}\label{eq:us}
\bu=\bJ\bs +v_0\left[\begin{array}{c}\bone\\ \bzero\end{array}\right]\quad\text{where} \quad 
\bJ:=\left[\begin{array}{cc}\bR & \bX\\ \bX & -\bR\end{array}\right]
\end{equation}
with $(\bR,\bX)$ as defined in \eqref{eq:LDF}. Note that \eqref{eq:LDF} is in fact the top block of \eqref{eq:us}. We therefore obtain $\nabla_{\bs}\bu=\bJ$ at $\bu=\bbu$.
\end{enumerate}
Plugging \emph{c1)} and \emph{c2)} in \eqref{eq:Hij} provides $\bH_n=\bJ^\top\tbN_n\bJ$.

We have hitherto derived the quadratic approximation $L_n^{(2)}(\bs) = \bs^\top\bJ^\top\tbN_n\bJ\bs$ for the losses on line $n$. Ohmic losses across all lines can be now approximated as
\begin{equation}\label{eq:L2}
L^{(2)}(\bs)=\sum_{n}L_n^{(2)}(\bs)=\frac{1}{2}\bs^\top\bJ^\top\tbN\bJ\bs
\end{equation}
where matrix $\tbN$ is defined as
\begin{equation*}
\tbN:=\sum_{n}\tbN_n=\left[\begin{array}{cc}
\bN & \bzero\\
\bzero & \bN
\end{array}\right]\quad \text{and}\quad \bN:=\sum_{n}2g_n\bgamma_n\bgamma_n^\top.
\end{equation*}
A key point here is that $\bN$ can be simplified as
\begin{equation}\label{eq:bN}
\bN=2\bGamma^\top\bD_r(\bD_r^2+\bD_x^2)^{-1}\bGamma.
\end{equation}
Thanks to \eqref{eq:bN}, we will next show that the matrix appearing in \eqref{eq:L2} features the simple form
\begin{equation*}
\bJ^\top\tbN\bJ=\left[\begin{array}{cc}
2\bR & \bzero\\
\bzero & 2\bR
\end{array}\right].
\end{equation*}
Once this form is plugged into \eqref{eq:L2}, it proves Proposition~\ref{pro:losses}.

To prove the above structure, we commence with the two diagonal blocks of $\bJ^\top\tbN\bJ$. Using the block partitions of $\bJ$ and $\tbN$, these two blocks are both equal to
\begin{align}\label{eq:db}
\bR\bN\bR+\bX\bN\bX&=\left(\bGamma^{-1}\bD_r\bGamma^{-\top}\right)\bN\left(\bGamma^{-1}\bD_r\bGamma^{-\top}\right)\nonumber\\
&~~~+\left(\bGamma^{-1}\bD_x\bGamma^{-\top}\right)\bN\left(\bGamma^{-1}\bD_x\bGamma^{-\top}\right)\nonumber\\
&=2\bGamma^{-1}\bD_r(\bD_r^2+\bD_x^2)^{-1}\bD_r^2\bGamma^{-\top}\nonumber\\
&~~~+2\bGamma^{-1}\bD_r(\bD_r^2+\bD_x^2)^{-1}\bD_x^2\bGamma^{-\top}\nonumber\\
&=2\bGamma^{-1}\bD_r(\bD_r^2+\bD_x^2)^{-1}(\bD_r^2+\bD_x^2)\bGamma^{-\top}\nonumber\\
&=2\bGamma^{-1}\bD_r\bGamma^{-\top}\nonumber\\
&=2\bR.
\end{align}

The top-right block of $\bJ^\top\tbN\bJ$ equals $\bR\bN\bX-\bX\bN\bR$. Using matrix manipulations similar to those in \eqref{eq:db}, it can be easily shown that this block is zero. The bottom-left block of $\bJ^\top\tbN\bJ$ equals $\bX\bN\bR-\bR\bN\bX$, so it is zero as well.
\color{black}

\section{Proof of Proposition~\ref{prop:exact penalty}}\label{sec:appendixB}
Let $\bg(\bx)$ collect $g_i(\bx)$ for all $i$. An optimal primal-dual pair of \eqref{eq:P1} satisfies \cmb{$\bg(\bx)\leq\bzero$, $\blambda\geq\bzero$,} and the conditions for Lagrangian optimality and complementary slackness:
\begin{subequations}
\begin{align}
&\hbx = \arg\min_{\bx\in\mcX}f(\bx)+\hblambda^\top\bg(\bx)\label{con:P1-1}\\
&\hblambda^\top\bg(\hbx)=\bzero.\label{con:P1-2}
\end{align}
\end{subequations}

Let $\tblambda$ denote the vector of Lagrange multipliers for \eqref{eq:P2}. Using the KKT conditions, we will show that $(\tbx,\ts;\tblambda)=(\hbx,0;\hblambda)$ is an optimal primal-dual pair for~\eqref{eq:P2}. Since $\bg(\hbx)\leq \bzero$, $\hbx\in\mcX$, and $\hblambda\geq\bzero$, the suggested triplet $(\hbx,0;\hblambda)$ is primal-dual feasible for \eqref{eq:P2}. Complementary slackness follows trivially from \eqref{con:P1-2}. Lagrangian optimality for \eqref{eq:P2} requires that $(\tbx,\ts;\tblambda)$ should satisfy
	\begin{align}
	\left(\hbx,\ts\right)= \arg\min_{\bx\in\mcX,s\geq 0}f(\bx)+\tblambda^\top\bg(\bx)-s\cdot \tblambda^\top \bone+p(s).\nonumber
	\end{align}	
Recognizing this problem is separable and for $\tblambda=\hblambda$ yield
\begin{subequations}
		\begin{align}
	&\tbx=\arg\min_{\bx\in\mcX} f(\bx)+\hblambda^\top\bg(\bx)\label{eq:P2:a}\\
	&\ts=\arg\min_{s\geq 0}~p(s)-s\cdot \hblambda^\top\bone.\label{eq:P2:b}
	\end{align}
\end{subequations}
Using \eqref{con:P1-1}, condition \eqref{eq:P2:a} provides that $\tbx=\hbx$. 
Because $\frac{d p(s)}{ds}> \hblambda^\top\bone$, the cost in \eqref{eq:P2:b} is increasing, 
and so $\ts=0$.

\balance
\bibliographystyle{IEEEtran}
\bibliography{myabrv,power}
\end{document}